\documentclass[11pt]{amsart}
\usepackage[T1]{fontenc}
\usepackage{lmodern,amsmath,amssymb,amsthm}
\usepackage[a4paper,margin=28mm]{geometry}
\usepackage{xcolor,hyperref}
\definecolor{LinkBlue}{HTML}{002FA7}
\definecolor{CitationRed}{HTML}{800020}
\definecolor{UrlGray}{HTML}{304C68}
\hypersetup{colorlinks=true,linkcolor=LinkBlue,
  citecolor=CitationRed,urlcolor=UrlGray}
\numberwithin{equation}{section}
\newtheorem{theorem}{Theorem}[section]
\newtheorem{proposition}[theorem]{Proposition}
\newtheorem{lemma}[theorem]{Lemma}
\newtheorem{corollary}[theorem]{Corollary}
\theoremstyle{remark}

\newcommand{\X}{\mathcal X}
\newcommand{\V}{\mathcal V}
\newcommand{\Z}{\mathcal Z}
\newcommand{\Y}{\mathcal Y}
\newcommand{\F}{\mathcal F}
\newcommand{\G}{\mathcal G}
\newcommand{\cO}{\mathcal O}
\newcommand{\Q}{\mathbb Q}
\newcommand{\C}{\mathbb C}
\newcommand{\db}{\bar\partial}
\newcommand{\ceil}[1]{\lceil #1\rceil}
\DeclareMathOperator{\Div}{div}
\DeclareMathOperator{\Supp}{Supp}
\DeclareMathOperator{\Fix}{Fix}
\hypersetup{pdftitle={Invariance of log plurigenera in families of compact Kähler manifolds},
  pdfauthor={Zhengyu Hu, Kefeng Liu}}
\title[Invariance of log plurigenera]{Invariance of log plurigenera in families of compact K\"ahler manifolds}
\author{Zhengyu Hu}
\address{Mathematical Sciences Research Center, Chongqing University of Technology,
Chongqing 400054, China}
\email{zhengyuhu16@gmail.com}
\author{Kefeng Liu}
\address{Mathematical Sciences Research Center, Chongqing University of Technology,
Chongqing 400054, China}
\email{kefengliu@cqut.edu.cn}
\date{November 26, 2026}
\begin{document}
\begin{abstract}
Let $(X_t,B_t)$ be a smooth proper family of pairs over a connected complex manifold with K\"ahler fibres, where the boundary has
relative simple normal crossings and real coefficients in
$[0,1)$. We show that for every $m\ge1$, every section of
$mK_{X_t}+\lfloor mB_t\rfloor$ extends locally over the base. In particular, the
log plurigenera are locally constant, and for $B_t=0$ we obtain the deformation
invariance of the plurigenera $P_m(X_t)$.
Since the total space need not be K\"ahler, we run the minimal model program on the central fibre, using the work of Hacon--Wang and Hacon--Xie, and extend its steps over a disc.
Vanishing along the fixed part on a common resolution gives an obstruction representative satisfying the integrability hypothesis of Cao--P\u{a}un.
As an application, if the adjoint is pseudoeffective on one fibre, then the pair has a
relative good minimal model near each point of the base and a relative
canonical model over the whole base.
\end{abstract}
\maketitle

\section{Introduction}
Siu conjectured that the plurigenera $P_m(X_t)=h^0(X_t,mK_{X_t})$ of a smooth family of compact K\"ahler manifolds over a disc are independent of $t$ \cite[Conjecture~0.4]{Siu}. We prove this conjecture, and the same statement for log plurigenera.
Let $p\colon\X\to S$ be a proper surjective holomorphic submersion
of complex manifolds with connected fibres $X_t=p^{-1}(t)$, $t\in S$. We say that a divisor has
\emph{relative simple normal crossings} if its components are smooth over
$S$ and, locally, their union is a union of fibre coordinate hyperplanes.
Boundaries have real coefficients unless we say otherwise.
For a boundary $\mathcal B$ on $\X$, we write $B_t=\mathcal B|_{X_t}$.

\begin{theorem}\label{main}
Suppose that $S$ is connected, each $X_t$ is K\"ahler, and $\mathcal B$ has
relative simple normal crossings with coefficients in $[0,1)$. Then for every integer $m\ge1$,
$$
                   h^0(X_t,mK_{X_t}+\lfloor mB_t\rfloor)
$$
is independent of $t$. The direct image
$p_*\cO_{\X}(mK_{\X/S}+\lfloor m\mathcal B\rfloor)$ is locally free and commutes
with holomorphic base change. Every such section on $X_t$ extends over
$p^{-1}(U)$ for each Stein open neighbourhood $U$ of $t$.
\end{theorem}

Taking $\mathcal B=0$, we obtain the following corollary.

\begin{corollary}\label{plurigenera}
Let $S$ be connected and $p\colon\X\to S$ a smooth proper family of
connected compact K\"ahler manifolds. Then for every $m\ge1$, the plurigenus
$P_m(X_t)=h^0(X_t,mK_{X_t})$ is independent of $t\in S$.
The local freeness, base change and extension over Stein neighbourhoods
in Theorem~\ref{main} apply to $\omega_{\X/S}^{\otimes m}$.
\end{corollary}

If the total space is K\"ahler, the
invariance was proved by Li--P\u{a}un--Wang--Zhou \cite[Theorem~1.2]{LPWZ}, and the logarithmic
statement in every rounded degree by Liu--Wan \cite[Theorem~1.1 and Corollary~1.2]{LW}.
But a family with K\"ahler fibres does not necessarily have K\"ahler total
space, even locally over the base: see Deligne's example \cite[Example~3.9]{Bingener}, recalled in \cite[\S1]{RT}. So in Theorem~\ref{main} the hypothesis of these two results is not satisfied in general.

For compact
complex surfaces, Iitaka proved invariance under every holomorphic deformation \cite[Theorem~III]{Iit}.
Levine \cite{Lev83,Lev85} studied infinitesimal extension of pluricanonical forms, and
Nakayama \cite{Nak} showed that, in the algebraic setting, invariance follows from the relative minimal model
conjectures.
For smooth projective families, Siu proved invariance first in general type
\cite{Siu98} and then in all cases \cite[Corollary~0.3]{Siu}.
Rao--Tsai treated smooth families of Moishezon manifolds
\cite[Main Theorem~1.2(i)]{RT}, and Hacon--Li--Rao proved the threefold case
for families with K\"ahler total space \cite[Theorem~6.7]{HLR}.

For the logarithmic version in log smooth projective families, Hacon--McKernan--Xu proved invariance in degrees clearing the boundary when
the adjoint or the boundary is relatively big
\cite[Theorem~1.8(1)]{HMX13}. Their general theorem \cite[Theorem~4.2]{HMX}
relies on the extension theorem of Berndtsson--P\u{a}un \cite[Theorem~0.2]{BP}.
Closer to our approach, Feng--Li--Rao deformed a central minimal model program. They assume that the central
pair is projective and canonical with big adjoint, and that its boundary
has no common component with the negative part \cite[Theorem~1.6]{FLR}.

\smallskip
\noindent\textbf{Minimal models.}\enspace We write $\mathcal D=K_{\X/S}+\mathcal B$
and $D_t=K_{X_t}+B_t$. By the \emph{relative canonical model} we mean the ample
model of $\mathcal D$ in the sense of \cite[Definitions~3.6.5 and~3.6.7]{BCHM},
taken locally on the analytic base. Note that a good minimal model maps to it by
the contraction defined by its semiample adjoint.

\begin{theorem}[Local good minimal model and relative canonical model]\label{model-theorem}
In Theorem~\ref{main}, suppose $D_s$ is pseudoeffective for some $s\in S$.
Then the following hold.
\begin{enumerate}
\item Each point of $S$ has an analytic neighbourhood $U$ over which
$(\X,\mathcal B)$ admits a good minimal model. More precisely, there is
a $\mathcal D$-negative bimeromorphic contraction
$\phi\colon\X_U\dashrightarrow\Y$ with $\Y$ strongly $\Q$-factorial,
$(\Y,\phi_*\mathcal B)$ klt and $K_{\Y/U}+\phi_*\mathcal B$ semiample over $U$.
\item There is a relative canonical model $\X\dashrightarrow\Z$ over $S$.
For every $t\in S$ its fibre is the canonical model of $(X_t,B_t)$.
If $\mathcal B$ is rational, the relative algebra
$\bigoplus_{m\geq0}p_*\cO_{\X}(\lfloor m\mathcal D\rfloor)$ is locally
finitely generated and $\Z$ is its relative $\operatorname{Proj}$.
\end{enumerate}
\end{theorem}

Hacon--Wang give termination and a good minimal model \cite[Theorem~1.1 and Corollary~1.2]{HW}.
For semiampleness they use the K\"ahler abundance theorem cited there as [OAI3, Theorem~1.2],
with subadditivity supplied by OAI4 or Hacon--Popa--Schnell. We use Hacon--Xie for contractions
and rational singularities \cite[Theorem~1.4 and Lemma~2.41]{HX}, Fujino for relative
$\operatorname{Proj}$ \cite[Theorem~1.8]{Fuj}, and adapt the construction of \cite[\S3]{HLX}
for strong $\Q$-factorialization near a compact fibre.

A relative good minimal model over all of $S$ does not exist in general.
In Proposition~\ref{global-counterexample} we construct a locally trivial family of
smooth projective threefolds over an elliptic curve $S$ and let
$\mathcal B=\tfrac12\mathcal E$. On each fibre we have $|2jD_t|=\{jE_t\}$ for
$j\geq1$, where $E_t=\mathcal E|_{X_t}$, so the relative canonical model
is $S$ itself. A $\mathcal D$-nonpositive bimeromorphic contraction with
relatively semiample transformed adjoint must contract $\mathcal E$.
But the monodromy is a nontorsion translation, and because of it every line bundle
on $\mathcal E$ has degree zero on a general fibre of the elliptic
surface $E_t$. So $\mathcal E$ cannot be contracted over $S$.
On the other hand, if $S$ is compact and
$\mathcal D$ is relatively big, then a global good minimal model exists
(Corollary~\ref{gmm-compact-big}).

\smallskip
\noindent\textbf{Sketch of proof.}\enspace We work over a disc with parameter
$t$ and consider a nonzero log pluricanonical section of degree $m\ge2$ on
the central fibre. In a projective family one
can run the relative MMP with scaling. Once the adjoint is relatively semiample, its free part on a common resolution then gives the positivity for
Koll\'ar-type injectivity in $H^1$, as in the proof of \cite[Lemma~6.3]{HMX}.
If only the fibres are K\"ahler, the total space need not be K\"ahler, and we have no relative MMP for it. Instead we run the program on the central fibre
and extend its steps over a disc. In Section~\ref{models} we round the
boundary in degree $m$, separate its components by blowups, and subtract
the fixed boundary contribution. By Lemma~\ref{reduction}, we only need to
extend sections of $mD$ on the new central pair $(X,B)$, where
$D=K_X+B$. This pair is canonical and $N_\sigma(D)\wedge B=0$. Since the program contracts only components of $N_\sigma(D)$, no
boundary component is contracted. By Hacon--Wang
\cite[Theorem~1.1 and Corollary~1.2]{HW}, there is a finite program ending in
a pair $(Y,B_Y)$ with $D_Y=K_Y+B_Y$ semiample. To extend the steps we use the deformation
theorem of Hacon--Xie \cite[Theorem~7.7]{HX}, which goes back to Koll\'ar
\cite[Theorem~2]{Kol} and, for a projective centre, to Hacon--Li--Rao \cite[Proposition~4.16]{HLR}. In this way we obtain relative contractions
and flips whose central fibres are the steps we started with (Proposition~\ref{relmodels}).
Next we look at a common resolution $X\xleftarrow{\mu}W\xrightarrow{\nu}Y$. We let
$F=\mu^*D-\nu^*D_Y$ and denote by $J$ the Jacobian divisor of $\mu$.
We pick $q>m$ sufficiently divisible that $|qD_Y|$
is free. A general member yields a section $\tau$ of $qD$ with
$\Div(\mu^*\tau)=qF+H$, where $qF$ is the fixed part of $|q\mu^*D|$
and $H$ meets the strata of $\Supp(J+F+\mu^*B)$ transversely.
For $k=q,m$, the obstruction to lifting a degree-$k$ jet from modulo $t^N$ to modulo $t^{N+1}$ lies in $H^1(X,kD)$.
On successive models the jets are the same
(Lemma~\ref{jets}), and so are the images of their obstructions on $W$.
Since $D_Y$ is semiample, $\nu^*D_Y$ has no fixed part. Comparing with $Y$, we find a representative vanishing along $J+\ceil{kF}$
(Proposition~\ref{obstruction}).
The key observation is that the same fixed part appears in the weight for lifting the jet. Indeed, on $W$ we have
$$
 \mu^*\bigl((k-1)\Div(\tau)/q+B\bigr)=(k-1)F+\tfrac{k-1}{q}H+\mu^*B .
$$
Since the representative vanishes along $J+\ceil{kF}$, $(X,B)$ is klt and
$(k-1)/q<1$ along $H$, the obstruction representative
on $X$ is square-integrable (Lemma~\ref{integrable}). So by the infinitesimal extension
theorem of Cao--P\u{a}un \cite[Theorem~1.2]{CP}, the jet lifts by one order.
Repeating this, we lift to every order, and by coherence and Krull intersection we get a holomorphic extension of the central section. Note that we first extend $\tau$ with its own
weight, and only then lift the given section of degree $m$ with the help of the extension of $\tau$.
For Theorem~\ref{model-theorem} we cannot use the reduction of Lemma~\ref{reduction}: subtracting fixed components of the boundary can change the minimal model. So we keep the original rational boundary and lift the central program step by step. If a lifted step has a different central fibre, some discrepancy goes up, and we restart on the new central pair. Only finitely many valuations have negative initial discrepancy, and their new discrepancies at restarts lie in $\frac1d\mathbb Z\cap(-1,0]$, where $d\mathcal B$ is
integral, so there are only finitely many restarts; by Theorem~\ref{main}, sections lift at every stage (Proposition~\ref{gmm-local}). For real coefficients, we approximate the boundary by the vertices of a rational simplex on which the asymptotic fixed orders are affine (Lemma~\ref{adjoint-box}). Their moving parts descend to a common model, and a rational interior boundary identifies the canonical models on each overlap.

\smallskip
\noindent\textbf{AI disclosure.}
The authors took the lead in developing the mathematical ideas and proof
strategies, with assistance from AI tools in working out parts of the proofs.
AI tools drafted the text, which the authors revised and reviewed in full. The authors
have verified all results and take full responsibility for the manuscript.

\section{Extending the central MMP}\label{models}
We work with complex analytic spaces as in \cite[\S\S2--3]{Fuj}.
A variety is reduced and irreducible. K\"ahler forms have smooth strictly
plurisubharmonic potentials in local embeddings \cite[Definition~2.2]{DHY}.
For generalized pairs on normal analytic spaces we use
\cite[Definition~2.1]{HLR} and \cite[Definition~2.7]{DHY}.

In the notation of Theorem~\ref{main}, fix $m\ge1$ and restrict to a disc
$\Delta$ in $S$ with coordinate $t$. Put $X=p^{-1}(0)$, $B=B_0=\mathcal B|_X$
and $n=\dim X$. We shrink $\Delta$ about $0$ as needed. For a model $\V$
over $\Delta$, write $V$ for its central fibre. A section modulo $t^N$
is a section on the closed subspace $t^N=0$.

We use additive notation for line bundles. On a normal space $V$,
$\cO_V(L)$ denotes the rank-one reflexive sheaf associated to an
integral Weil divisor $L$. For an effective integral Cartier divisor
$G$, we write $\sigma_G$ for its canonical section.

Replacing $\mathcal B$ by $\lfloor m\mathcal B\rfloor/m$, we may assume that
$m\mathcal B$ is integral. If $H^0(X,m(K_X+B))=0$, upper semicontinuity suffices.
For $m=1$ the rounded boundary is zero and $h^{n,0}$ is locally constant
by Hodge decomposition and upper semicontinuity. So we can assume $m\ge2$
and $|m(K_X+B)|\ne\varnothing$.

\begin{lemma}\label{reduction}
After a sequence of blowups of relative strata $\rho\colon\widetilde{\X}\to\X$,
there is a boundary $\mathcal B'$ on $\widetilde{\X}$ with coefficients in
$[0,1)\cap\frac1m\mathbb Z$ whose components are disjoint and smooth
over $\Delta$. On the central fibre $\widetilde X$ put $B'=\mathcal B'|_{\widetilde X}$
and $D'=K_{\widetilde X}+B'$. No component of $B'$ is fixed in $|mD'|$.
In particular, $N_\sigma(D')\wedge B'=0$. There is a natural isomorphism
$$
 H^0\big(\widetilde X,mD'\big)
 \simeq H^0\big(X,m(K_X+B)\big).
$$
If a section on the left extends over $\Delta$, so does its image.
\end{lemma}
\begin{proof}
We blow up the intersection of two boundary components of coefficients $a,b$
and give the exceptional divisor $R$ coefficient $\max(0,a+b-1)$. The new
adjoint equals the pullback of the old one plus $\max(0,1-a-b)R$, so pushforward
preserves degree-$m$ sections on the family and on each fibre.
The fibres remain K\"ahler since they are blowups along smooth centres.
Each pair is separated once, and each new positive coefficient is
strictly smaller than both $a$ and $b$. By descending induction on the
coefficients, the process terminates.

Let $\mathcal B^*$ be the resulting boundary, $B^*=\mathcal B^*|_{\widetilde X}$ and $F=\Fix|m(K_{\widetilde X}+B^*)|$.
We extend $E=mB^*\wedge F$ along the relative components to $\mathcal E$ and put
$\mathcal B'=\mathcal B^*-\mathcal E/m$. Multiplication by $\sigma_E$ identifies the central
sections, and multiplication by $\sigma_{\mathcal E}$ carries their extensions back.
The new fixed part $F-E$ shares no component with $B'$. Since
$N_\sigma(D')\le \frac1m\Fix|mD'|=(F-E)/m$, we obtain $N_\sigma(D')\wedge B'=0$.
\end{proof}

We henceforth write $(\X,\mathcal B)$ for the pair of Lemma~\ref{reduction} and
put $\mathcal D=K_{\X/\Delta}+\mathcal B$, $D=K_X+B$. Since the boundary
components are disjoint, $(X,B)$ is canonical.

The smooth compact K\"ahler pair $(X,B)$ is klt with rational boundary,
and $D$ has an effective multiple. Scaling a K\"ahler class so that its
sum with $D$ is nef, we obtain by \cite[Theorem~1.1 and the proof of
Corollary~1.2]{HW} a finite minimal model program
$$
 (X,B)=(V_0,B_0)\dashrightarrow\cdots\dashrightarrow
 (V_l,B_l)=(Y,B_Y),\qquad D_i=K_{V_i}+B_i ,
$$
with $D_Y=K_Y+B_Y$ semiample.

\begin{proposition}\label{relmodels}
Each central step $(V,B)\xrightarrow{h_0}Z\xleftarrow{h_0^+}(V^+,B^+)$
extends to a diagram of normal varieties flat over $\Delta$,
\begin{equation}\label{diagram}
        (\V,\mathcal B)\xrightarrow{h}\Z\xleftarrow{h^+}(\V^+,\mathcal B^+) .
\end{equation}
Both arrows are projective, and the fibre at $0$ is the given diagram.
For a divisorial contraction the right arrows are identities and
$\mathcal B^+=h_*\mathcal B$. For a flip both arrows are small and
$\mathcal B^+$ is the strict transform of $\mathcal B$.
In either case $\mathcal D^+=K_{\V^+/\Delta}+\mathcal B^+$ is
$\Q$-Cartier, $(\V^+,V^++\mathcal B^+)$ is plt, and
$$
 (K_{\V^+}+V^++\mathcal B^+)|_{V^+}=K_{V^+}+B^+ .
$$
\end{proposition}
\begin{proof}
By induction, $(\V,V+\mathcal B)$ is plt with different $B$ on $V$.
The initial case follows from log smoothness. Since $-(K_V+B)$ is $h_0$-ample,
$$
 (h_0)_*\cO_V=\cO_Z,\qquad R^1(h_0)_*\cO_V=0 .
$$
Pulling back the diagram of \cite[Proposition~11.4.1]{KM} along a classifying
map for $\V\to\Delta$ extends $h_0$ to a contraction
$h\colon\V\to\Z$ with $\Z$ flat over $\Delta$ and $\Z\times_\Delta\{0\}=Z$.
Its target is normal since $\V$ is normal and $h_*\cO_\V=\cO_\Z$.
The fibres over $U=Z\setminus h_0(\operatorname{Exc}h_0)$ are points.
By upper semicontinuity of fibre dimension, $h$ is finite over a neighbourhood
of $U$ in $\Z$, hence an isomorphism there.

Since $V$ is compact, the adjoint $\mathcal D=K_{\V/\Delta}+\mathcal B$
has a Cartier multiple $r\mathcal D$ near $V$.
By adjunction, $(-r\mathcal D)|_V=-r(K_V+B)$ is $h_0$-ample. By the analytic
fibrewise ampleness criterion \cite[\S1.1, Theorems~3.2.7 and~3.2.9]{Con},
$-r\mathcal D$ is $h$-ample near $Z$, hence over a smaller disc by properness.
So $h$ is projective and bimeromorphic, so $\mathcal D$ is
$h$-big. Over small Stein open sets, \cite[Theorem~1.8(2),(3)]{Fuj}
gives the relative log canonical model
$$
 \V^+=\operatorname{Proj}_{\Z}\bigoplus_{d\ge0}
                   h_*\cO_{\V}(dr\mathcal D).
$$
The adjoint $\mathcal D^+$ is $\Q$-Cartier and $h^+$-ample, so a common
Cartier multiple makes $h^+$ projective as above. The normal space $\V^+$ dominates
the disc and hence is flat. The model extracts no divisor and is an isomorphism
near $U$. Thus $T=\Div_{\V^+}(t)$ is a reduced prime Cartier divisor.
Take a common log resolution of $(\V,V+\mathcal B)$ and
$(\V^+,T+\mathcal B^+)$ with projections $\alpha$ and $\beta$. By $\mathcal D$-nonpositivity
\cite[Lemma~11.3(3)]{Fuj},
\begin{equation}\label{step-pullback}
    \alpha^*\mathcal D=\beta^*\mathcal D^++E,\qquad E\ge0,\quad \beta_*E=0 .
\end{equation}
Since $\alpha^*V=\beta^*T=\Div(t)$, $(\V^+,T+\mathcal B^+)$ is plt.
Hence $T$ is normal and
$K_T+\Gamma=(K_{\V^+}+T+\mathcal B^+)|_T$ with $\Gamma\ge0$.

Let $\widetilde V$ be the strict transform of $V$, with induced maps
$a\colon\widetilde V\to V$ and $b\colon\widetilde V\to T$.
By adjunction on $\widetilde V$, we have
\begin{equation}\label{central-pullback}
 a^*(K_V+B)=b^*(K_T+\Gamma)+E|_{\widetilde V},\qquad E|_{\widetilde V}\ge0 .
\end{equation}
The central MMP contracts only components of $N_\sigma(D_0)$
\cite[Lemma~5.6]{HX}. Since $N_\sigma(D_0)\wedge B_0=0$, every boundary
component survives, so monotonicity of discrepancies carries the canonicity
of $(V_0,B_0)$ to $(V,B)$. Thus $a^*(K_V+B)\le K_{\widetilde V}+\widetilde B$,
where $\widetilde B$ is the strict transform of $B$.
Every component of $B$ meets $h_0^{-1}(U)$, where the boundaries agree.
Pushing the discrepancy inequality forward by $b$, we have
$$
 \Gamma+b_*(E|_{\widetilde V})\le b_*\widetilde B\le\Gamma .
$$
Hence $\Gamma=b_*\widetilde B$ and the effective divisor $E|_{\widetilde V}$
is $b$-exceptional. Pushing down \eqref{central-pullback}, we therefore obtain,
for $u=h^+|_T$ and every sufficiently divisible $d>0$,
$$
 (h_0)_*\cO_V\big(d(K_V+B)\big)=u_*\cO_T\big(d(K_T+\Gamma)\big).
$$
Since $K_T+\Gamma$ is $u$-ample, relative $\operatorname{Proj}$ identifies
$T$ with the given ample model $V^+$ over $Z$. Under this identification
$b_*\widetilde B=B^+$.

In the divisorial case $h^+$ is finite, hence an isomorphism, near $Z$.
In a flip any exceptional divisor would be horizontal and cut out a
contracted divisor on the central fibre. Hence both arrows are small.
\end{proof}

\section{Jet obstructions and the fixed part}\label{obstructions}
Recall that $N_\sigma(D_0)\wedge B_0=0$ for the auxiliary pair
$(V_0,B_0)=(X,B)$, and that $D_Y$ is semiample on the final model.
Write $(\V_i,\mathcal B_i)\to\Delta$ for the families of
Proposition~\ref{relmodels}. By adjunction,
$\mathcal D_i=K_{\V_i/\Delta}+\mathcal B_i$ restricts to
$D_i=K_{V_i}+B_i$ on $V_i=(t=0)$.
On a common resolution $a\colon\mathcal W\to\X$, $b_i\colon\mathcal W\to\V_i$,
by successive applications of \eqref{step-pullback} we have
$a^*\mathcal D-b_i^*\mathcal D_i\ge0$ with $b_i$-exceptional support.

For $k\ge2$ with $k\mathcal B$ integral, this identifies the degree-$k$
sections over every open subset of $\Delta$. Put
$\F_i=\cO_{\V_i}(k\mathcal D_i)$ and $\F_{i,N}=\F_i/t^N\F_i$.
Koll\'ar's theorem \cite[Theorem~2 and (4.1)]{Kol11} applied to the integral
$\Q$-Cartier divisor $k\mathcal D_i$ on the plt pair of Proposition~\ref{relmodels}
shows that $\F_i$ is Cohen--Macaulay, so its quotient $\F_{i,1}$ by the
nonzerodivisor $t$ satisfies $S_2$. Since $\cO_{V_i}(kD_i)$ is reflexive on
the normal fibre $V_i$, the adjunction isomorphism in codimension one extends
uniquely to
\begin{equation}\label{restriction}
                     \F_{i,1}\simeq\cO_{V_i}(kD_i).
\end{equation}

\begin{lemma}\label{jets}
For a step \eqref{diagram}, write $\F,\F^+$ for the adjoint sheaves on
$\V,\V^+$ and $\F_N,\F_N^+$ for their quotients modulo $t^N$.
Put $\G=h_*\F$ on $\Z$. For every integer $N\ge1$,
$$
       h_*\F_N=\G/t^N\G=h^+_*\F^+_{N}.
$$
Taking global sections gives
\begin{equation}\label{jet-sections}
 H^0(\V,\F_N)\simeq H^0(\Z,\G/t^N\G)\simeq H^0(\V^+,\F_N^+).
\end{equation}
The identifications respect the quotient $\F_{N+1}\to\F_N$ and
multiplication $\F_N\xrightarrow{\,t\,}\F_{N+1}$.
\end{lemma}
\begin{proof}
Push forward the short exact sequence
$$
 0\longrightarrow\F\xrightarrow{\,t^N\,}\F
 \longrightarrow\F_N\longrightarrow0 .
$$
From the connecting map, we obtain a second short exact sequence
$$
 0\longrightarrow\G/t^N\G\longrightarrow h_*\F_N
 \longrightarrow (R^1h_*\F)[t^N]\longrightarrow0 .
$$
Here $[t^N]$ denotes the kernel of multiplication by $t^N$. So restriction
$h_*\F\to h_*\F_N$ is surjective exactly when this kernel vanishes.
The same criterion applies to $h^+$.

On the positive side, $R^1h^+_*\F^+=0$. This is immediate when $h^+$ is
the identity. For a flip, $h^+$ is projective, $(\V^+,\mathcal B^+)$ is klt,
and $(k-1)\mathcal D^+$ is $h^+$-ample. Hence Fujino's vanishing theorem
applies to the integral $\Q$-Cartier divisor $k\mathcal D^+$
\cite[Theorem~5.2]{Fuj}. In particular, restriction to the central fibre
is surjective on the positive side.

The differences in \eqref{step-pullback} and \eqref{central-pullback} are
effective and exceptional over the respective positive models. So
$h_*\F=h^+_*\F^+$ on $\Z$ and, using \eqref{restriction},
$h_*\F_1=h^+_*\F^+_1$ on $Z$.
Under these identifications both restriction maps reduce the same sections
modulo $t$. Hence, by the sequence for $N=1$, we have
$$
 (R^1h_*\F)[t]
 =\operatorname{coker}\bigl(\G/t\G\longrightarrow h^+_*\F^+_1\bigr)=0 .
$$
So every $t^N$ acts injectively on $R^1h_*\F$. By the second sequence,
we obtain the asserted identifications, with both compatibilities induced by restriction.
\end{proof}

On the central fibre $X=V_0$, resolve the maps $X\dashrightarrow V_i$
by successive blowups with smooth centres \cite[\S1.2.4]{AKMW}. We obtain
a compact K\"ahler manifold $W$ and morphisms $\pi_i\colon W\to V_i$.
Write $\mu=\pi_0$, $\nu=\pi_l$ and set
$$
 J=K_W-\mu^*K_X,\qquad E_i=\mu^*D-\pi_i^*D_i,\qquad F=E_l .
$$
By negativity, $E_i\ge0$ with $\pi_i$-exceptional support, and $J$ is
also $\pi_i$-exceptional since $X\dashrightarrow V_i$ extracts no divisor.
We take $W$ so that $\Supp(J+F+\mu^*B)$ has simple normal crossings.

Set $M_k=K_W+\mu^*(kD-K_X)=k\mu^*D+J$.
For each $i$ and $N\ge1$, multiplication by $t^N$ and \eqref{restriction} give
$$
 0\longrightarrow\cO_{V_i}(kD_i)\xrightarrow{\,t^N\,}\F_{i,N+1}
 \longrightarrow\F_{i,N}\longrightarrow0 ,
$$
since $t$ is a nonzerodivisor on $\F_i$. Write
$\delta_{i,N}\colon H^0(\V_i,\F_{i,N})\to H^1(V_i,kD_i)$
for its connecting homomorphism.
For a morphism $f\colon T'\to T$ and a coherent sheaf $\mathcal H$ on $T'$,
write $f^*$ for the edge map in the low-degree Leray sequence
\begin{equation}\label{leray-edge}
 0\longrightarrow H^1(T,f_*\mathcal H)\xrightarrow{\,f^*\,}H^1(T',\mathcal H)
 \longrightarrow H^0(T,R^1f_*\mathcal H).
\end{equation}
For $\mu$, we have $\mu_*\cO_W(M_k)=\cO_X(kD)$, and the injection
$H^1(X,kD)\hookrightarrow H^1(W,M_k)$ is pullback followed by multiplication by
$\sigma_J$, where $J=\Div(\det d\mu)$.

\begin{proposition}\label{obstruction}
Let $N\ge1$. Then
$$
 \operatorname{im}(\mu^*\circ\delta_{0,N})\ \subseteq\
 \operatorname{im}\!\left(H^1(W,M_k-J-\ceil{kF})\longrightarrow H^1(W,M_k)\right).
$$
\end{proposition}
\begin{proof}
We show that the obstruction on $W$ comes from $Y$, where the fixed part has been removed. Let $T_k=M_k-J-\ceil{kF}$. Then we have
$$
       \pi_{i*}\cO_W(M_k)=\cO_{V_i}(kD_i),\qquad
       \nu_*\cO_W(T_k)=\cO_Y(kD_Y).
$$
Indeed, since $M_k=\pi_i^*(kD_i)+kE_i+J$ and $T_k=\lfloor\nu^*(kD_Y)\rfloor$,
the local sections in each equality are identified by pullback.

Let $s_{0,N}\in H^0(\X,\F_{0,N})$. At each step, by \eqref{jet-sections}
we send $s_{i,N}$ to $u_{i,N}\in H^0(\Z,\G/t^N\G)$ and then to $s_{i+1,N}$.
We write $\xi_i=\delta_{i,N}(s_{i,N})$ and $\G_r=\G/t^r\G$.
By Lemma~\ref{jets}, both lifting sequences push forward to
$0\to\G_1\xrightarrow{\,t^N\,}\G_{N+1}\to\G_N\to0$.
Let $\zeta\in H^1(Z,\G_1)$ be the connecting image of $u_{i,N}$.

We next compare $\zeta$ with $\xi_i$ and $\xi_{i+1}$. To this end, let $f\colon V_j\to Z$ be either map of the central step, where $j=i$ or $i+1$.
We view the jet sheaves as abelian sheaves on the central fibres.
The identification $\G_r\simeq f_*\F_{j,r}$ is adjoint to a map
$e_r\colon f^{-1}\G_r\to\F_{j,r}$. At $x$ it sends a section over
$f^{-1}(U)$, for $U$ a neighbourhood of $f(x)$, to its germ at $x$.
Since $(f^{-1}\G_r)_x=(\G_r)_{f(x)}$ for each $x\in V_j$, the upper row below is exact.
$$
\begin{array}{ccccccccc}
0&\to& f^{-1}\G_1 &\xrightarrow{\,t^N\,}& f^{-1}\G_{N+1}
 &\longrightarrow& f^{-1}\G_N &\to&0\\
 &&\big\downarrow{\scriptstyle e_1}&&\big\downarrow{\scriptstyle e_{N+1}}
 &&\big\downarrow{\scriptstyle e_N}\\
0&\to& \F_{j,1} &\xrightarrow{\,t^N\,}& \F_{j,N+1}
 &\longrightarrow& \F_{j,N} &\to&0 .
\end{array}
$$
For either arrow $a\colon\F_{j,r}\to\F_{j,s}$ in the bottom row, both paths
around its square are adjoint to $f_*a\colon\G_r\to\G_s$. Hence the diagram
commutes. By
\eqref{jet-sections}, pullback followed by $e_N$ sends $u_{i,N}$ to $s_{j,N}$. Since $e_1$ is adjoint to the identification
$\G_1=f_*\F_{j,1}$, the induced map in degree one is the natural comparison
map for $\Gamma(Z,-)\circ f_*$, hence the Leray edge map \eqref{leray-edge}
by \cite[Theorem~5.8.3 and \S5.8.6]{Wei}. By naturality of the connecting maps,
$f^*\zeta=\delta_{j,N}(s_{j,N})=\xi_j$.
Since the two maps $W\to Z$ agree,
$$
 \pi_i^*\xi_i=(h_0\circ\pi_i)^*\zeta=(h_0^+\circ\pi_{i+1})^*\zeta
             =\pi_{i+1}^*\xi_{i+1}
 \qquad\text{in }H^1(W,M_k).
$$
Finally, we pass to $Y$. Note that the inclusion $\cO_W(T_k)\hookrightarrow\cO_W(M_k)$ pushes forward under
$\nu$ to the identity on $\cO_Y(kD_Y)$. So by naturality of the edge map in
\eqref{leray-edge}, the image of $\xi_l$ in $H^1(W,T_k)$ maps to
$\nu^*\xi_l=\mu^*\delta_{0,N}(s_{0,N})$ in $H^1(W,M_k)$. This completes the proof.
\end{proof}

\section{Extension of sections}\label{extension}
Choose $q>m$ sufficiently divisible that $q\mathcal B$ and $qF$ are integral and
$qD_Y$ is a globally generated line bundle. For a general section $\tau_Y$
of $qD_Y$, let $\mu^*\tau=\sigma_{qF}\nu^*\tau_Y$ and $H=\nu^*\Div(\tau_Y)$.
Then $\Div(\mu^*\tau)=qF+H$. Applying analytic
Bertini \cite[Theorem~3.2]{Fuj23} to $\nu^*|qD_Y|$ on $W$ and on the strata of $\Supp(J+F+\mu^*B)$
shows that $H$ is smooth (possibly empty) and transverse to every stratum.
We first extend $\tau$ in degree $q$, then use its extension to extend an
arbitrary section in degree $m$. For either degree $k$ put
$$
 \Theta_k=\frac{k-1}{q}\Div(\tau)+B,\qquad
 \mu^*\Theta_k=(k-1)F+\frac{k-1}{q}H+\mu^*B .
$$
So $\mu$ is a log resolution of $(X,\Theta_k)$. We define the multiplier ideal by
$$
 \mathcal J(\Theta_k):=\mu_*\cO_W(J-\lfloor\mu^*\Theta_k\rfloor).
$$
This ideal is independent of the chosen log resolution
\cite[(5.8)--(5.9)]{Dem}.

\begin{lemma}\label{integrable}
Let $N\ge1$ and let $s_N$ be a holomorphic section of $k\mathcal D$
modulo $t^N$. Then there exist a smooth section $w$ and a smooth
$(0,1)$-form $\lambda$, both with values in $k\mathcal D$ near $X$, such that
$$
 w\equiv s_N\pmod{t^N},\qquad \db w=t^N\lambda .
$$
The restriction $\lambda_X=\lambda|_X$ represents the obstruction
$\delta_{0,N}(s_N)$ to lifting $s_N$ modulo $t^{N+1}$. Near each point
of $X$ it has a finite expression
$\lambda_X=\sum_\alpha f_\alpha\beta_\alpha$, where
$f_\alpha\in\mathcal J(\Theta_k)$ and $\beta_\alpha$ are smooth
$kD$-valued $(0,1)$-forms.
\end{lemma}
\begin{proof}
On each open set $U_i$ of a cover near $X$, choose a holomorphic section
$u_i$ representing $s_N$ modulo $t^N$. On overlaps $u_j-u_i$ is divisible
by $t^N$. For a subordinate partition of unity $\{\rho_j\}$, set
$w_0=\sum_j\rho_j u_j$. On $U_i$ we have
$$
  w_0=u_i+t^N v_i,\qquad
  v_i=\sum_j\rho_j\,\frac{u_j-u_i}{t^N}.
$$
Since $u_i$ is holomorphic, $\db w_0=t^N\db v_i$ on $U_i$. Hence
$\lambda_0=t^{-N}\db w_0$ extends smoothly across $X$.

The map $a|_X\mapsto t^N a\bmod t^{N+1}$ identifies $\cO_X(kD)$
with the kernel of $\F_{0,N+1}\to\F_{0,N}$. With $\F_0=\cO_\X(k\mathcal D)$,
the lifting sequence is the bottom row of the following diagram near $X$.
$$
\begin{array}{ccccccccc}
0&\to&\F_0&\xrightarrow{\,t^N\,}&\F_0&\to&\F_{0,N}&\to&0\\
&&\big\downarrow{\scriptstyle\mathrm{res}_X}
&&\big\downarrow{\scriptstyle\bmod\,t^{N+1}}&&\big\Vert\\
0&\to&\cO_X(kD)&\xrightarrow{\,t^N\,}&\F_{0,N+1}&\to&\F_{0,N}&\to&0 .
\end{array}
$$
In the top row, apply $\db$ to the smooth lift $w_0$ of $s_N$.
Its differential $\db w_0=t^N\lambda_0$ is the image of $\lambda_0$
under multiplication by $t^N$. The form $\lambda_0$ is $\db$-closed since
$t^N\db\lambda_0=\db^2w_0=0$, so the upper connecting map sends $s_N$
to $[\lambda_0]$. The left vertical arrow restricts this class to $X$,
giving the connecting class in the bottom row,
$$
 \delta_{0,N}(s_N)=[\lambda_0|_X].
$$
Changing the smooth lift by $t^Nh$ changes the representative
$\lambda_0|_X$ by the exact form $\db(h|_X)$.

All coefficients of $\mu^*B-J$
are less than one by klt, and $\mu^*\Theta_k$ has coefficient
$(k-1)/q<1$ along $H$. Hence $\lfloor\mu^*\Theta_k\rfloor\le J+\ceil{kF}$.
Since $\mu$ is bimeromorphic, $-\mu^*\Theta_k$ is $\mu$-nef and $\mu$-big.
By \cite[Corollary~1.4]{Fuj13} and the projection formula, we have
$R^j\mu_*\cO_W(J-\lfloor\mu^*\Theta_k\rfloor)=0$ for $j>0$.
By Proposition~\ref{obstruction} and the Leray spectral sequence,
$\delta_{0,N}(s_N)$ lifts to
$$
 H^1(W,M_k-\lfloor\mu^*\Theta_k\rfloor)\simeq H^1\big(X,\cO_X(kD)\otimes\mathcal J(\Theta_k)\big).
$$

We pick an acyclic cover $\{U_i\}$ of $X$ and a subordinate partition
of unity $\{\rho_j\}$. We represent the lifted class by a \v{C}ech cocycle
$c_{ij}\in\Gamma(U_i\cap U_j,\cO_X(kD)\otimes\mathcal J(\Theta_k))$,
with $c_{ij}+c_{j\ell}=c_{i\ell}$. On $U_i$ set
$$
 \lambda_X|_{U_i}=\sum_j\db\rho_j\,c_{ij}.
$$
The cocycle identity shows that these forms agree on overlaps. Indeed,
the local smooth sections $b_i=\sum_j\rho_jc_{ij}$ satisfy
$b_i-b_j=c_{ij}$ and $\db b_i=\lambda_X|_{U_i}$, so the
\v{C}ech--Dolbeault comparison identifies $[\lambda_X]$ with the image
of the lifted class in $H^1(X,\cO_X(kD))$. Since each $c_{ij}$ is a
section of $\cO_X(kD)\otimes\mathcal J(\Theta_k)$, the displayed formula
yields the stated local finite-sum expression for $\lambda_X$.

It remains to choose a smooth lift $w$ of $s_N$ for which
$(t^{-N}\db w)|_X=\lambda_X$. Since $[\lambda_0|_X]=[\lambda_X]$, write
$\lambda_0|_X-\lambda_X=\db v$ for a smooth section $v$ of $kD$ on $X$.
We extend $v$ smoothly to $\widetilde v$ near $X$ and put
$w=w_0-t^N\widetilde v$. Then $w\equiv s_N\pmod{t^N}$ and, with
$\lambda=\lambda_0-\db\widetilde v$, we have
$$
 \db w=t^N\lambda,\qquad \lambda|_X=\lambda_0|_X-\db v=\lambda_X .
$$
\end{proof}

\begin{theorem}[Cao--P\u{a}un {\cite[Theorem~1.2]{CP}}, special case]\label{CP-lifting}
Let $p\colon\X\to\Delta$ be smooth and proper with K\"ahler central
fibre $X$. For line bundles $\mathcal L,\mathcal E_1,\ldots,\mathcal E_s$
on $\X$, put $\mathcal A=K_{\X/\Delta}+\mathcal L$. Suppose
$$
 (1-r_0)\mathcal L=r_0K_{\X/\Delta}+\sum_{i=1}^s r_i\mathcal E_i
 \quad\text{in }\operatorname{Pic}(\X)\otimes\Q
 \qquad (r_0,\ldots,r_s\in\Q_{\ge0},\quad r_0<1).
$$
For $N\ge2$, let $u$ be a section of $\mathcal A$ modulo $t^N$ and
$v_i$ sections of $\mathcal E_i$ modulo $t^{N+1}$. Put $v_0=u$ and
$I=\{i\in\{0,\ldots,s\}:r_i>0\}$. Assume $v_i|_X\ne0$ for $i\in I$.
$$
 \Theta=\sum_{i\in I}r_i\Div(v_i|_X).
$$
Let $w\equiv u\pmod{t^N}$ be a
smooth lift with $\db w=t^N\lambda$, and put $\lambda_X=\lambda|_X$.
Suppose that near each point of $X$,
$\lambda_X=\sum_\alpha f_\alpha\beta_\alpha$ with
$f_\alpha\in\mathcal J(\Theta)$ and $\beta_\alpha$ smooth
$\mathcal A|_X$-valued $(0,1)$-forms. Then $u$ extends to a holomorphic
section modulo $t^{N+1}$.
\end{theorem}

We apply Theorem~\ref{CP-lifting} twice, first in degree $q$ and then in degree $m$. In both cases we take, for degree $k$,
$$
 \mathcal A=k\mathcal D,\qquad
 \mathcal L=k\mathcal D-K_{\X/\Delta},\qquad
 \Theta=\Theta_k .
$$
In degree $q$ we lift the current $q\mathcal D$-jet $u$ whose central value is $\tau$, and we let
$$
 r_0=\frac{q-1}{q},\qquad
 \mathcal E_1=q\mathcal B,\qquad
 r_1=\frac1q,\qquad
 v_1=\sigma_{q\mathcal B}.
$$
We see that $\Theta=(q-1)\Div(\tau)/q+B=\Theta_q$. After $\tau$ is extended, we lift the given
$m\mathcal D$-jet $u$ in degree $m$, and now we let
$$
 r_0=0,\qquad
 \mathcal E_1=q\mathcal D,\qquad
 \mathcal E_2=q\mathcal B,\qquad
 r_1=\frac{m-1}{q},\qquad r_2=\frac1q,\qquad
 (v_1,v_2)=(\widetilde\tau,\sigma_{q\mathcal B}).
$$
Here $\widetilde\tau$ denotes the holomorphic extension of $\tau$, and we have
$\Theta=(m-1)\Div(\tau)/q+B=\Theta_m$. Note that the relation
$$
 k\mathcal D-K_{\X/\Delta}
   =\frac{k-1}{q}(q\mathcal D)+\frac1q(q\mathcal B),
$$
is just $\mathcal L=r_0\mathcal A+\sum_i r_i\mathcal E_i$. So in both cases the
line bundle hypothesis of Theorem~\ref{CP-lifting} holds.

We also need the analytic hypotheses of \cite[Theorem~1.2]{CP} for the form $\lambda_X$ of
Lemma~\ref{integrable}. We write locally
$\lambda_X=\sum_\alpha f_\alpha\beta_\alpha$ as in that lemma. On a log
resolution $\pi\colon\widetilde X\to X$, by the definition of the multiplier ideal we have
$$
 \Div((f_\alpha\circ\pi)\det d\pi)\ge\lfloor\pi^*\Theta\rfloor .
$$
When we pull back the top-form component of $\lambda_X$ by $\pi$, the Jacobian
$\det d\pi$ appears. So by the inequality above,
$\pi^*\lambda_X$ is smoothly divisible by a local equation of
$\lfloor\pi^*\Theta\rfloor$. This is the condition \cite[(3.3.3)]{CP}.
We remove this integral part and change variables. What is left is
the fractional part of $\pi^*\Theta$, whose coefficients lie in $[0,1)$
since it has simple normal crossings support. Hence every weight
$|z|^{-2\{\gamma\}}$ we meet is locally integrable. We deduce that each
$f_\alpha\beta_\alpha$, and hence $\lambda_X$, is square-integrable with
weight $\Theta$. In particular, it is square-integrable with $(1-\varepsilon)\Theta$ for
$0<\varepsilon<1$, which is \cite[(1.2.1)]{CP}.

\begin{proof}[Proof of Theorem~\ref{main}]
We make the base change $t=z^2$ and write $\widetilde p\colon\widetilde\X\to
\Delta_z$ for the resulting family. Apply Sections~\ref{models}--\ref{obstructions} to this smooth family over $\Delta_z$ with the same central pair and central MMP. Its restriction to $z^2=0$ is the
product $X\times\operatorname{Spec}\C[z]/(z^2)$. The same holds for the
pulled-back adjoint bundles. So every central section has a constant
lift modulo $z^2$.
Starting with the chosen section $\tau$ of $qD$, we apply
Theorem~\ref{CP-lifting} successively to lift each jet modulo $z^N$ to one
modulo $z^{N+1}$. At every step the obstruction form satisfies the hypothesis of Theorem~\ref{CP-lifting}, by Lemma~\ref{integrable} and the discussion above.
Next we show that a central section which lifts to every order extends holomorphically. Let $a$ be a central section of degree $k\in\{q,m\}$, and we set
$\mathcal F=\cO_{\widetilde\X}(k\widetilde{\mathcal D})$ and
$M=(R^1\widetilde p_*\mathcal F)_0$. By Grauert's theorem, $M$ is a finite
$\C\{z\}$-module. We write $\mathcal F_N=\mathcal F/z^N\mathcal F$ and
$X_N=(z^N=0)$. The direct-image connecting maps for
$0\to\mathcal F\xrightarrow{z}\mathcal F\to\mathcal F_1\to0$ and
$0\to\mathcal F\xrightarrow{z^N}\mathcal F\to\mathcal F_N\to0$
give connecting classes $\delta(a),\delta_N(a_N)\in M$ for $a_N|_X=a$.
The map between these sequences is $z^{N-1}$ on the left terms, the
identity in the middle, and restriction on the quotients. Hence
$$
 \delta(a)=z^{N-1}\delta_N(a_N)\in z^{N-1}M .
$$
If $a$ lifts to every $X_N$, then $\delta(a)\in\bigcap_Nz^NM=0$ by
Krull's intersection theorem. Exactness of the first sequence gives a
holomorphic extension of $a$. Applying this to $\tau$, we obtain an extension
$\widetilde\tau$ of the chosen $qD$ section. Using $\widetilde\tau$ and
$\sigma_{q\mathcal B}$ in Theorem~\ref{CP-lifting}, we lift every jet of
any given $mD$ section. By the Krull argument above, this section also extends holomorphically.

Since the bundle is pulled back from the $t$-disc, we average an extension
with its pullback under $(x,z)\mapsto(x,-z)$ and obtain an invariant section.
It descends under $t=z^2$ with the same central value. Extending a basis, we see that $h^0(X_t,\cdot)\ge h^0(X_s,\cdot)$ near $s$.
By upper semicontinuity, equality holds. Restricting to complex-line discs,
we obtain local constancy on $S$. By Grauert's theorem, we have local freeness and
base change for $\mathcal G=p_*\cO_\X(mK_{\X/S}+\lfloor m\mathcal B\rfloor)$.
If $U$ is Stein, then
$H^1(U,\mathcal G\otimes\mathcal I_s)=0$ by Cartan's theorem B, hence
$H^0(U,\mathcal G)\to\mathcal G_s/\mathfrak m_s\mathcal G_s\simeq
H^0(X_s,mK_{X_s}+\lfloor mB_s\rfloor)$ is surjective. By openness of the
K\"ahler condition \cite{KS}, the conclusions hold locally near any K\"ahler fibre.
\end{proof}

\begin{corollary}\label{iitaka-volume}
In Theorem~\ref{main}, the invariant Iitaka dimension and the Iitaka volume
of $D_t=K_{X_t}+B_t$ are independent of $t$. The same holds for
pseudoeffectivity and the ordinary volume.
\end{corollary}
\begin{proof}
For rational $\mathcal B$, the growth of the plurigenera in
Theorem~\ref{main} determines the Iitaka dimension and both volumes, and
nonvanishing \cite[Corollary~1.2]{HW} detects pseudoeffectivity.
For real $\mathcal B$, we approximate from above by rational boundaries and obtain invariance of
pseudoeffectivity and of ordinary volume by closedness of the pseudoeffective
cone and continuity of volume \cite{Bou}.

For the Iitaka invariants we approximate effective representatives on fixed
supports. Let
$E\sim_{\mathbb R}D_t$ be an effective real representative and let $\kappa=\kappa(E)$, as in
\cite[Definitions~2.5--2.6]{Zhu}. Its Iitaka volume is
$\limsup_{n\to\infty}\kappa!\,h^0(X_t,\lfloor nE\rfloor)/n^\kappa$.
If there is no such $E$, the two invariants are $-\infty$ and $0$.
On a fixed support the Iitaka dimension is constant and the Iitaka volume
is continuous. Indeed, nearby effective divisors lie between
$(1-\varepsilon)E$ and $(1+\varepsilon)E$, and scaling $E$ by $c>0$
multiplies its Iitaka volume by $c^\kappa$.

We fix two fibres locally and pick an effective $E_0\sim_{\mathbb R}D_0$.
With the supports fixed, the relation $K_X+B_0\sim_{\mathbb R}E_0$
is rational affine in the coefficients, being a relation in
$\operatorname{Pic}(X)\otimes\mathbb R$. Hence $(\mathcal B,E_0)$ is a
convex combination of rational pairs $(\mathcal B_i,E_{0i})$ satisfying
the same relation and coefficient bounds, with weights $c_i>0$.
By the rational case, there are effective $E_{ti}\sim_{\Q}K_{X_t}+B_{i,t}$. For
positive rational weights $a_i$ of sum one, $\sum a_iE_{0i}$ and
$\sum a_iE_{ti}$ have equal Iitaka invariants by Theorem~\ref{main}.
Letting $a_i\to c_i$, we obtain equality for $E_0$ and
$\sum c_iE_{ti}\sim_{\mathbb R}D_t$ by fixed-support continuity.
If we impose two representative relations on the same fibre, the same argument shows that the invariants do not depend
on the choice of $E$. Exchanging the fibres, we obtain
the assertion when no effective representative exists.
\end{proof}

\section{Relative minimal and canonical models}\label{minimal-models}
We write $\mathcal D=K_{\X/S}+\mathcal B$ and $D_s=K_{X_s}+B_s$. A bimeromorphic contraction is
$\mathcal D$-negative if the difference of the pullbacks of the adjoints on a common resolution is
effective, exceptional over the target and positive on every prime lost
from the source. We use strong $\Q$-factoriality in the sense of \cite{HLX}:
each rank-one reflexive sheaf has an invertible reflexive power.

In this section we prove Theorem~\ref{model-theorem}: first for a rational boundary near one point (Proposition~\ref{gmm-local}), then for real boundaries, and then we glue. We use the following Proj comparison throughout. If
$a^*L=b^*L'+G$ on a common resolution, with $G\geq0$ $b$-exceptional,
then for divisible $k$, $b_*\cO_W(kG)=\cO_T$ and the projection formula
give $b_*\cO_W(ka^*L)=\cO_T(kL')$. Pushing to a common base, we identify
the corresponding relative section algebras and their ample models.

\subsection*{Central K\"ahler MMP}

\begin{lemma}\label{gmm-central}
Let $(V,B)$ be a compact K\"ahler klt pair with $B$ effective and rational
and $D=K_V+B$ pseudoeffective. Then there is a finite sequence ending in
a semiample adjoint, each step of which is a diagram
$$
 (V_i,B_i)\xrightarrow{h_i} Z_i \xleftarrow{h_i^+}(V_{i+1},B_{i+1})
$$
of projective bimeromorphic contractions of compact K\"ahler spaces.
Here $(V_0,B_0)=(V,B)$, $D_i=K_{V_i}+B_i$ is $\Q$-Cartier and the
pairs are klt. The divisor $-D_i$ is $h_i$-ample and the right side
is the relative ample model of $D_i$ over $Z_i$.
\end{lemma}

\begin{proof}
Since $V$ is not necessarily strongly $\Q$-factorial, we run the MMP of Hacon--Wang on a resolution $W$ of $V$ and show that the steps descend to models starting at $V$. Let $f\colon W\to V$ be a projective log resolution with an effective
exceptional divisor $C$ such that $-C$ is $f$-ample
\cite[proof of Theorem~3.3]{HLX}. We put
$B_W=f_*^{-1}B+(1-\epsilon)\sum P$ over the exceptional primes $P$.
For small rational $\epsilon>0$, since $(V,B)$ is klt, we have $D_W=K_W+B_W=f^*D+F$
where $F$ is positive on every exceptional prime. Let
$\omega$ be a K\"ahler class with $[D]+\omega$ K\"ahler. For $\delta>0$ small,
$\omega_W=f^*\omega-\delta[C]$ is K\"ahler and
\begin{equation}\label{gmm-scaling}
 [D_W]+\lambda\omega_W
 = f^*([D]+\lambda\omega)+[F-\lambda\delta C],\qquad
 F-\lambda\delta C\geq0\quad(0\leq\lambda\leq1).
\end{equation}
Near $\lambda=1$ this is the pullback of a K\"ahler class plus the class
of an effective $f$-exceptional divisor, so the ample model is $V$.

Next we run the MMP. The resolution $W$ is smooth and strongly $\Q$-factorial. Since $D_W$
is pseudoeffective, we can run the MMP of \cite[Theorem~1.1 and proof
of Corollary~1.2]{HW} on $(W,B_W)$ with scaling of a sufficiently
large multiple of $\omega_W$. It terminates with a semiample adjoint.
We use the parameter $\lambda$ of \eqref{gmm-scaling}. On a scaling
interval $[a,b]$, $a<b$, we write $P$ for the model, $L=K_P+B_P$ and
$\eta$ for the transform of $\omega_W$. Then
$\alpha_\lambda=[L]+\lambda\eta$ is nef on $[a,b]$. At a preceding
step with threshold $c\geq\lambda$, let $E\geq0$ be the difference of the pullbacks of the ordinary adjoints. The scaled difference is zero at $c$ and
affine in $\lambda$, hence equals $(1-\lambda/c)E$
\cite[\S3]{HW}. So the pair with nef part
$\lambda\overline{\omega_W}$ is generalized klt. For $\lambda>0$ its
boundary with nef part is modified big and $\alpha_\lambda$ is big by
\eqref{gmm-scaling} and birational pushforward. By
\cite[Theorem~1.4]{HX}, we obtain its contraction to a compact
K\"ahler space. At zero we use the semiample adjoint. At a common endpoint
the pullbacks agree across all steps at that value, so the endpoint
contractions agree \cite[Lemma~2.33]{HX}.

On each scaling interval we now construct one model $M$. Let $g\colon P\to M$ be the Stein factor of the product of the two
endpoint contractions. On $a<\lambda<b$ we have
$$
 \alpha_\lambda=\frac{b-\lambda}{b-a}\alpha_a+
                 \frac{\lambda-a}{b-a}\alpha_b=g^*\theta_\lambda ,
$$
where $\theta_\lambda$ is K\"ahler by restriction to the product image
and finite pullback. Hence $g$ is the bimeromorphic contraction throughout
this interval.

We show that $L$ descends to $M$. We pick two distinct interior values of $\lambda$. Their classes
$\alpha_\lambda$ are pullbacks, so extrapolating to $\lambda=0$ gives
$c_1(L)\in g^*H^2(M,\mathbb R)$. The image of $g^*$ is defined over
$\Q$, and $c_1(L)$ is rational. Clearing denominators and then torsion, we obtain a Cartier multiple $kL$ and a class $\xi\in H^2(M,\mathbb Z)$
with $c_1(kL)=g^*\xi$ in $H^2(P,\mathbb Z)$.

The map $g$ is proper and bimeromorphic between normal compact spaces in
Fujiki's class $\mathcal C$. The klt pair $(P,B_P)$ is generalized klt
with zero nef part, and $L=K_P+B_P$ has degree zero on every curve
contracted by $g$. By \cite[Lemma~2.41 and its proof]{HX},
$M$ has rational singularities. Since $P$ also has rational singularities,
$R^1g_*\cO_P=0$. On a contractible Stein neighbourhood $O\subset M$,
by Leray and Stein vanishing we have $H^1(g^{-1}O,\cO_P)=0$. Also
$c_1(kL)|_{g^{-1}O}=g^*(\xi|_O)=0$. By the exponential sequence, we see that
$\cO_P(kL)|_{g^{-1}O}\simeq\cO_{g^{-1}O}$. Its transition functions
descend since $g_*\cO_P=\cO_M$, so $kL$ is the pullback of a line bundle.
Pushing down, we identify it on $M$ with $\cO_M(kD_M)$, where
$D_M=g_*L=K_M+g_*B_P$. Thus $L=g^*D_M$ for compatible canonical
divisors. In particular, $(M,g_*B_P)$ is klt.

It remains to see that consecutive interval models are related by the steps in the statement. At a common endpoint $\tau>0$, let $M,M^+$ denote the interval models for
$\lambda>\tau$ and $\lambda<\tau$, respectively. Both map to the endpoint
space $Z_\tau$. On either side, by the identity $\alpha_\tau\equiv_{Z_\tau}0$ we have
$$
 [D_M]\equiv_{Z_\tau}-\frac{\tau}{\lambda-\tau}\theta_\lambda .
$$
Hence $-D_M$ and $D_{M^+}$ are ample over $Z_\tau$, and the maps
$h\colon M\to Z_\tau$ and $h^+\colon M^+\to Z_\tau$ are projective by
\cite[Lemma~2.39]{HX}. On a common resolution $u,v$, comparing the crepant
descents, we have
$$
 u^*D_M=v^*D_{M^+}+E_\tau,\qquad E_\tau\geq0,\quad v_*E_\tau=0 ,
$$
where $E_\tau$ is the sum of the pulled-back differences of the MMP steps at $\tau$.
Let $r>0$ clear the indices. By the Proj comparison, the
graded adjoint algebras over $Z_\tau$ are identified. Since $D_{M^+}$ is $h^+$-ample,
$M^+$ is their relative $\operatorname{Proj}$. Over a Stein
neighbourhood of each point of $Z_\tau$, the pair on $M$ is klt and
$D_M$ is relatively big because $h$ is bimeromorphic.
By \cite[Theorem~1.8(2),(3) and Definition~11.4(vi)]{Fuj}
the induced map $M\dashrightarrow M^+$ is the bimeromorphic contraction
to its relative log canonical model. It therefore extracts no divisor.

Finally, we check that the first model is $V$. For $\lambda$ near $1$, we compare \eqref{gmm-scaling} with the first
interval model on a common resolution. If the two decompositions are
$r_i^*\theta_i+[N_i]$, $N_i\geq0$ $r_i$-exceptional, then $N_1-N_2$
is $r_1$-nef with nonpositive pushforward. By negativity in both directions,
$N_1=N_2$, hence the same contraction
\cite[Lemmas~2.5 and~2.33]{HX}. The sequence starts at $V$ and ends with
$D_M$ semiample by $L=g^*D_M$ on the last interval. This completes the proof.
\end{proof}

\subsection*{Lifting one central contraction}

\begin{lemma}\label{gmm-flag}
Let $\phi\colon\X_U\dashrightarrow\V$ be a bimeromorphic contraction over a
neighbourhood $U$ of $s$, isomorphic near a dense open of $X_s$, with
$\V$ normal and proper over $U$. Suppose $\mathcal D_\V=K_{\V/U}+\phi_*\mathcal B$
is $\Q$-Cartier and the difference of the adjoint pullbacks is effective.
Then the central fibre $V$ is normal, $\V$ is flat near $V$, and adjunction gives
a klt pair $(V,B)$ with $K_V+B=\mathcal D_\V|_V$.

A diagram $V\to Z\leftarrow V^+$ as in Lemma~\ref{gmm-central} lifts,
after shrinking, to projective contractions $\V\to\Z\leftarrow\V^+$
with the given left side on the central fibre. The total step is
$\mathcal D_\V$-negative. Its positive central pair $(T,\Gamma)$
is klt and compact K\"ahler. On a common resolution $a\colon W\to V$, $b\colon W\to T$, we have
\begin{equation}\label{gmm-central-error}
 a^*(K_V+B)=b^*(K_T+\Gamma)+G,\qquad G\geq0 .
\end{equation}
If $b_*G=0$ in that formula, then $(T,\Gamma)$ is the given
positive pair over $Z$.
\end{lemma}
\begin{proof}
We first prove that $V$ is normal and klt and that $\V$ is flat, by cutting $\V$ down to $V$ with one coordinate hyperplane at a time. Let $t_1,\ldots,t_e$ be coordinates centred at $s$. We cut successively
by their zero divisors. Suppose that $D,D'$ are the adjoints of one pair of slices
and that $a^*D=b^*D'+G'$ with $G'\geq0$ on a common log resolution. For the
next cut, we denote by $H,T'$ the source and target divisors of $t_j$. Then
$$
 a^*(D+H)=b^*(D'+T')+G',\qquad a^*H=b^*T'=\Div(t_j).
$$
Note that the source with $H$ added is plt and $H$ is its only lc place. Comparing discrepancies via the equality, we see that the target is also plt, and its only lc
place is the transform of $H$. Since the two sides agree near the general point of $H$,
$T'$ is reduced and prime. By adjunction, it is normal with klt different.
Restricting the equality on the resolution to the strict transform of $H$, we see that the difference of the pullbacks is still effective. So by induction all slices along the flag are
normal and klt, and all differences are effective. Moreover, the coordinate
divisors form a regular sequence. Hence by the local flatness criterion
\cite[cf.\ Proposition~8]{KNX}, the family is flat over $U$.

Next we lift the step. For the left contraction $h_0\colon V\to Z$ we have $R^1(h_0)_*\cO_V=0$.
So by
\cite[Proposition~11.4.1]{KM}, we can argue as in the proof of
Proposition~\ref{relmodels}: the extension and the relative Proj construction work over $U$. After shrinking, we obtain a projective
$h\colon\V\to\Z$ with central fibre $h_0$ and $-r\mathcal D_\V$
$h$-ample. Its positive side is
$$
 h^+\colon\V^+=\operatorname{Proj}_{\Z}
             \bigoplus_{j\geq0}h_*\cO_\V(jr\mathcal D_\V)\longrightarrow\Z
$$
by \cite[Theorem~1.8(2),(3)]{Fuj}. As in that proposition, we see that the
map $\V\dashrightarrow\V^+$ is a bimeromorphic contraction with effective pullback error. For a prime lost from $\V$, a general curve on its strict transform contracted over $\Z$ moves through its general point and has curve image in $\V$. The negative-side adjoint has negative degree on this curve and the positive-side adjoint has nonnegative degree. The error thus contains the strict transform of the prime, giving $\mathcal D_\V$-negativity.

It remains to look at the central fibre of $\V^+$. We apply the flag argument to the composite from $\X_U$. This gives the normal slices and their adjunction boundaries on $\V^+$. On a common log resolution of the current step $\V\dashrightarrow\V^+$, restrict its effective adjoint error successively to the strict transforms of the same coordinate slices. Their general points lie in the isomorphism locus away from the error. Each restriction is therefore effective, and the pullbacks of the coordinate divisor agree at every cut. Adjunction gives \eqref{gmm-central-error} for this step.
Note that $T$ is projective over $Z$, so it is compact K\"ahler. Now assume that $b_*G=0$.
Then for sufficiently divisible $k$,
the relative adjoint algebras are the same by the Proj comparison. Hence
$T\simeq V^+$, and by this formula, $\Gamma=B^+$ for compatible
canonical divisors.
\end{proof}

\subsection*{Strong $\Q$-factorialization}

\begin{lemma}\label{gmm-factorial}
Let $\Y\to U$ be proper over a coordinate ball and let $(\Y,\Theta)$
be klt. After shrinking about its centre there is a projective small modification
$\Y'\to\Y$ with $\Y'$ strongly $\Q$-factorial.
\end{lemma}

\begin{proof}
We first make finitely many rank-one reflexive sheaves $\Q$-Cartier, and then show that finitely many suffice, because $H^2$ of a small tube around the central fibre is finite dimensional. Any finitely many rank-one reflexive sheaves on $\Y$ can be
made $\Q$-Cartier by a projective small modification near the central
fibre. To see this, we use the construction of \cite[Proposition~3.1]{HLX} for one
sheaf at a time. For a sheaf $L$, the local small modifications are
the relative ample models of its reflexive algebra
$\bigoplus_{m\geq0} L^{[m]}$, hence they glue on intersections. Since the central fibre is compact, after shrinking there is a common Cartier multiple and the modification is relatively projective. Note that later small modifications keep
the $\Q$-Cartier transforms we already have.

Now we explain why finitely many sheaves are enough. Let $f:R\to\Y$ be a projective resolution and let $g:R\to U$ be the
composite. We claim that for sufficiently small concentric balls $U'$, the groups
$H^2(R_{U'},\Q)$ are all the same finite-dimensional space under restriction.
Indeed, by properness the critical values of $\rho=\sum_j|t_j\circ g|^2$ near $0$
are subanalytic \cite[\S3]{BM}, and by
Sard's theorem they have measure zero. Hence some interval $(0,\epsilon]$ consists of regular
values. Between its levels we have product collars, and so the open tubes have the homotopy
type of a compact sublevel manifold.

In this stable cohomology group, the first Chern classes of line bundles on all such tubes span a rational subspace, and we pick a basis of it. We represent the basis
by $L_1,\ldots,L_r$ on one tube. Then we apply the first paragraph to
$(f_*L_j)^{**}$, and we denote their $\Q$-Cartier transforms on
$\mu\colon\Y'\to\Y$ by $A'_j$.

We show that $\Y'$ is strongly $\Q$-factorial. Let $A'$ be a rank-one reflexive sheaf on $\Y'$. We put $A=(\mu_*A')^{**}$
and $L=(f^*A)^{**}$. Clearing the denominators of $c_1(L)$ in the chosen
basis, the cohomological torsion and the indices of the $A'_j$, we see that there are
integers $n>0,n_j$ and a line bundle $P$ on $\Y$ with
$$
 c_1\!\left(L^n\otimes\bigotimes_j L_j^{-n_j}\right)=0,\qquad
 L^n\simeq f^*P\otimes\bigotimes_j L_j^{n_j}.
$$
For the second equality we use the exponential sequence; note that
$H^1(\Y,\cO_\Y)\simeq H^1(R,\cO_R)$ because $\Y$ has rational singularities.
Now we restrict to the complement of the exceptional loci and extend reflexively across
codimension two on $\Y'$. We obtain
$(A')^{[n]}\simeq\mu^*P\otimes\bigotimes_j(A'_j)^{[n_j]}$, which is a line bundle. This completes the proof.
\end{proof}

\subsection*{Local good minimal models for rational boundaries}

\begin{lemma}[Bound on the number of restarts]\label{gmm-restarts}
For a sequence of steps lifted from $(X_s,B_s)$ as in Lemma~\ref{gmm-flag}, choose $d>0$ such that
$d\mathcal B|_{\X_U}$ is integral near $X_s$. Let $\mathcal P_s$ be the
set of divisorial valuations $P$ over $X_s$ with $a(P,X_s,B_s)<0$.
Then $\mathcal P_s$ is finite, and the number of lifted steps with
$b_*G\ne0$ in \eqref{gmm-central-error} is at most
$$
                 \sum_{P\in\mathcal P_s} \bigl(-d\,a(P,X_s,B_s)\bigr).
$$
\end{lemma}
\begin{proof}
Let $f\colon\widetilde X\to X_s$ be the separating blowups in
Lemma~\ref{reduction}, and take the boundary $B^*$ there before
subtracting the fixed part. Comparing the pullbacks of the adjoints, we have
$$
 K_{\widetilde X}+B^*=f^*(K_{X_s}+B_s)+E,\qquad E\geq0 .
$$
The components of $B^*$ are smooth and disjoint, hence $(\widetilde X,B^*)$
is canonical. For $P\in\mathcal P_s$, we have
$a(P,\widetilde X,B^*)=a(P,X_s,B_s)-\operatorname{ord}_P E<0$.
So $P$ is a component of $B^*$ and $\mathcal P_s$ is finite.

By \eqref{gmm-central-error}, along every lifted step the discrepancy
of each valuation is nondecreasing. For a prime component $P$ of $b_*G$
at a restart, the increase is
\begin{equation}\label{gmm-discrepancy}
 a(P,T,\Gamma)-a(P,V,B)=\operatorname{coeff}_P(b_*G)>0 .
\end{equation}
The normal fibre $T$ is smooth at a general point of $P$, and flatness
makes the family smooth there. With $b_j$ the coefficients of the
surviving total boundary components $\mathcal B_j^+$, by adjunction we have
\begin{equation}\label{gmm-coefficients}
 -a(P,T,\Gamma)=\operatorname{coeff}_P\Gamma
 =\sum_j b_j\operatorname{mult}_P(\mathcal B_j^+|_T)
 \in\tfrac1d\mathbb Z\cap[0,1).
\end{equation}
Consequently $a(P,X_s,B_s)\leq a(P,V,B)<a(P,T,\Gamma)\leq0$.
For a fixed $P$, the new discrepancies at the restarts assigned to $P$
are distinct points of $(a(P,X_s,B_s),0]\cap\tfrac1d\mathbb Z$.
The initial discrepancy lies in the same lattice since $dD_s$ is Cartier.
There are at most $-d\,a(P,X_s,B_s)$ such values. Summing over $P$, we obtain the bound.
\end{proof}

\begin{proposition}\label{gmm-local}
Suppose in Theorem~\ref{main} that $\mathcal B$ is rational and $D_s$
is pseudoeffective. Then on an analytic neighbourhood $U$ of $s$ there is a
$\mathcal D$-negative bimeromorphic contraction $\phi\colon\X_U\dashrightarrow\Y$
with $\Y$ strongly $\Q$-factorial, $(\Y,\phi_*\mathcal B)$ klt and
$K_{\Y/U}+\phi_*\mathcal B$ semiample over $U$.
\end{proposition}
\begin{proof}
The idea is to run the central MMP, lift it step by step, and restart whenever the lifted step has a different central fibre. We work over a Stein neighbourhood of $s$ and pick a finite central program
as in Lemma~\ref{gmm-central}. For a step starting at $(V,B)$,
Lemma~\ref{gmm-flag} constructs a relative step with central pair
$(T,\Gamma)$. Let $G$ be the effective divisor in \eqref{gmm-central-error}.
There are two cases. If $b_*G=0$, the Proj comparison identifies the adjoint algebras over
$Z$ and hence their relative $\operatorname{Proj}$. So $T$ is the given
positive side, and we go on to the next step. If $b_*G\ne0$, the central fibre has changed, and we restart the program from $(T,\Gamma)$.

For any finite chain, let $(\Y,\mathcal B_\Y)$ be the current total pair,
with central fibre $T$, and write $D_\Y=K_{\Y/U}+\mathcal B_\Y$ and
$D_T=D_\Y|_T$. On a common resolution $a_s\colon W_s\to X_s$,
$b_s\colon W_s\to T$, the differences add up to
$G_s=a_s^*D_s-b_s^*D_T\geq0$. For divisible $q$ clearing the indices,
negativity gives $H^0(\Y,qD_\Y)\simeq H^0(\X_U,q\mathcal D)$, while the
central map is the injection
$$
 j_q\colon H^0(T,qD_T)\longrightarrow H^0(X_s,qD_s),\qquad
 v\longmapsto\sigma_{qG_s}b_s^*v .
$$
These maps commute with restriction on the common dense open. By Theorem~\ref{main},
restriction from $\X_U$ onto $X_s$ is surjective. Hence $j_q$ is an
isomorphism, and restriction from $\Y$ onto $T$ is surjective.

Now $H^0(X_s,qD_s)\ne0$ for divisible $q$ by \cite[Corollary~1.2]{HW}.
Since $j_q$ is an isomorphism, every restart again has pseudoeffective adjoint. By Lemma~\ref{gmm-central},
we obtain a finite program each time, and by Lemma~\ref{gmm-restarts} the
number of restarts is bounded. So the process stops, and on the final fibre $D_T$ is semiample. We pick a free multiple of $D_T$.
By surjectivity, its basis lifts. By properness, the evaluation is surjective
after shrinking. Finally, Lemma~\ref{gmm-factorial} provides a small crepant
modification. It preserves negativity and semiampleness.
\end{proof}

\subsection*{Real boundaries and relative canonical models}

Write $\mathcal B=\sum_{i=1}^h b_i\mathcal B_i$ with $b_i>0$. For a vector
$c$ put $\mathcal B(c)=\sum_i c_i\mathcal B_i$, $\mathcal D(c)=K_{\X/S}+\mathcal B(c)$,
and denote their restrictions by $B_t(c),D_t(c)$. We fix rational endpoints
$0<b_i^-<b_i<b_i^+<1$ and let $q$ be twice a common denominator.

\begin{lemma}\label{adjoint-box}
Each $s\in S$ has an analytic neighbourhood on which the multigraded algebra
\begin{equation}\label{box-algebra}
 \bigoplus_{\substack{d\in q\mathbb Z_{\geq0},\ \alpha\in\mathbb Z^h\\
                      0\leq\alpha_i\leq d(b_i^+-b_i^-)}}
 p_*\cO_{\X}\Bigl(d\mathcal D(b^-)+\sum_i\alpha_i\mathcal B_i\Bigr)
\end{equation}
is finitely generated and its fibre at $t$ is the corresponding section
ring on $X_t$. The locus of $c$ in the box $\prod_i[b_i^-,b_i^+]$ with
$D_t(c)$ pseudoeffective is independent of $t$. When nonempty it is a
rational polytope with a subdivision into rational cells on which all
asymptotic orders of the rational systems vary affinely in $c$, both on
the family and on every fibre.
\end{lemma}
\begin{proof}
Let $\pi\colon P\to\X$ be the fibre product of the quotient bundles
$\mathbf P_{\X}(\cO\oplus\cO(\mathcal B_i))$. We write $H_i$ for the
$i$-th tautological class on $P$. In that factor the quotient killing
$\cO(\mathcal B_i)$ gives a section $T_i^\infty$ of class
$H_i-\pi^*\mathcal B_i$. The vectors $(1,0)$ and
$(1,\sigma_{\mathcal B_i})$ cut out two sections $T_i^0,T_i^1$ of class
$H_i$, disjoint from $T_i^\infty$. Over a local equation $z_i=0$ of
$\mathcal B_i$ the latter two have equations $w_i=0$ and $w_i+z_i=0$.
Since $\X$ has relative normal crossings, so do all these divisors on $P$. We set
$$
 \Gamma=\sum_i\left((1-b_i^-)T_i^\infty+
                   \tfrac{1+b_i^+}{2}(T_i^0+T_i^1)\right).
$$
Each $P_t$ is compact K\"ahler and $(P_t,\Gamma_t)$ is klt. The line bundle
$L=q(K_{P/S}+\Gamma)$ satisfies
$$
 L\sim q\pi^*\mathcal D(b^-)+\sum_i q(b_i^+-b_i^-)H_i .
$$
For $d=kq$, by the projective bundle formula we have
$$
 \pi_*\cO_P(kL)=
 \bigoplus_{0\leq\alpha_i\leq d(b_i^+-b_i^-)}
 \cO_\X\Bigl(d\mathcal D(b^-)+\sum_i\alpha_i\mathcal B_i\Bigr).
$$
This decomposition commutes with restriction to the fibres, and the product
of the summands indexed by $\alpha,\beta$ has index $\alpha+\beta$.
Applying $p_*$, we obtain \eqref{box-algebra}. Let
$r=p\circ\pi\colon P\to S$ and let $\C(s)$ be the residue field at $s$.
By Theorem~\ref{main} for $(P,\Gamma)$, base change holds in every positive
degree, and $r_*\cO_P=\cO_S$ since the fibres are connected.
Since restriction respects products, we obtain an isomorphism of graded algebras
$$
 \C(s)\otimes_{\cO_{S,s}}
 \left(\bigoplus_{k\geq0}(r_*\cO_P(kL))_s\right)
 \ \xrightarrow{\ \sim\ }\ \bigoplus_{k\geq0}H^0(P_s,kL_s).
$$
The right side is a Veronese of the finitely generated adjoint ring of
$(P_s,\Gamma_s)$ by \cite[Theorem~1.8]{Fuj15}. Hence the proper map $r$
of normal spaces, the line bundle $L$ and the closed subset $\{s\}$
satisfy \cite[Lemma~11]{Kol}. By this lemma, there is one neighbourhood of $s$ on
which the whole graded algebra is finitely generated.

We pick local frames $\sigma_i$ in the finitely many positive generator
degrees. They generate on the family and restrict to nonzero generators
on every fibre. We write $w_i$ for their weights, $C=\sum_i\mathbb R_{\geq0}w_i$
and $w(c)=(1,c-b^-)$. For rational $c$, some multiple of $w(c)$ is
the weight of a monomial exactly when $w(c)\in C$. By nonvanishing
\cite[Corollary~1.2]{HW}, this is equivalent to pseudoeffectivity of $D_t(c)$.
Both loci are closed. In the pseudoeffective locus rational points
approximate each $c$ from above within the box, and in the degree-one
slice of $C$ rational points are dense by polyhedrality. Hence the
equivalence holds for real $c$ as well.

For a divisorial valuation $P$ and a rational weight $w\in C$, the
asymptotic order is
\begin{equation}\label{weight-order}
 o_P(w)=\min_{\substack{\lambda_i\geq0\\ \sum_i\lambda_iw_i=w}}
                          \sum_i\lambda_i\operatorname{ord}_P(\sigma_i).
\end{equation}
Indeed, monomials span each weight space, and clearing the denominators
of a rational minimizer realizes this order in a multiple of $w$.
By \cite[Lemma~2.1]{ELMNP} the weights $w_i$ determine a finite rational
fan on $C$ on which all these minima are linear. The same formula on
$X_t$ uses $\sigma_i|_{X_t}$ with the same weights. Slicing the fan by
$w=w(c)$, we obtain the asserted subdivision.
\end{proof}

\begin{proof}[Proof of Theorem~\ref{model-theorem}(1)]
We apply Proposition~\ref{gmm-local} at the barycentre of a small rational simplex around $b$, and show that this model is good for every point of its relative interior. Let $s\in S$. By Corollary~\ref{iitaka-volume}, $D_s(b)$ is pseudoeffective.
Let $A$ be the smallest rational affine subspace containing $b$. Let $C_b$
be the unique cell of the rational subdivision in Lemma~\ref{adjoint-box}
whose relative interior contains $b$. Its affine span is rational and
contains $A$. Hence $C_b\cap A$ contains a neighbourhood of $b$ relative
to $A$, and we can choose a simplex there with rational vertices
$v_1,\ldots,v_r$ and $b$ in its interior relative to $A$.
We apply Proposition~\ref{gmm-local} to its barycentre $b_*=r^{-1}\sum_jv_j$.
We denote by $\phi\colon\X_U\dashrightarrow\Y$ the model so obtained and put
$\mathcal D_{\Y}(c)=K_{\Y/U}+\phi_*\mathcal B(c)$. On a projective common
resolution $a\colon W\to\X_U$, $g\colon W\to\Y$
resolve also the base ideals of sufficiently divisible multiples on the
vertex rays. Each ray algebra is finitely generated by Gordan's lemma
applied to the monomials of that ray. So we can write
$$
 a^*\mathcal D(v_j)=M_j+F_j,\qquad F_j\geq0,
$$
where $M_j$ is relatively semiample and $F_j$ is the asymptotic fixed
divisor. At $b_*$ semiampleness identifies the fixed divisor $F_*$ with
the difference of the pullbacks. Since the fixed divisors vary affinely on the simplex, we have
\begin{equation}\label{common-fixed}
 F_*:=a^*\mathcal D(b_*)-g^*\mathcal D_{\Y}(b_*)=\tfrac1r\sum_j F_j,\qquad
 g^*\mathcal D_{\Y}(b_*)=\tfrac1r\sum_j M_j .
\end{equation}
$F_*$ is $g$-exceptional, so every $F_j$ is $g$-exceptional.
Since $\phi$ extracts no divisor, $g_*M_j=\mathcal D_{\Y}(v_j)$ in
compatible frames, and the latter is $\Q$-Cartier by strong
$\Q$-factoriality. Hence $M_j-g^*\mathcal D_{\Y}(v_j)$ is exceptional
and $g$-nef. By negativity and \eqref{common-fixed}, we have
$$
 M_j-g^*\mathcal D_{\Y}(v_j)\leq0,\qquad
 \sum_j\bigl(M_j-g^*\mathcal D_{\Y}(v_j)\bigr)=0 .
$$
Each difference vanishes. Relative generation descends by
$g_*\cO_W=\cO_\Y$, so every $\mathcal D_\Y(v_j)$ is semiample over $U$.
For an interior $c=\sum_j t_jv_j$ all $t_j>0$, and
$$
 a^*\mathcal D(c)-g^*\mathcal D_\Y(c)=\sum_j t_jF_j,\qquad
 \Supp\Bigl(\sum_j t_jF_j\Bigr)=\Supp\Bigl(\sum_j F_j\Bigr).
$$
So strict negativity holds as at $b_*$, and klt follows by discrepancy
comparison. The adjoint $\mathcal D_\Y(c)=\sum_jt_j\mathcal D_\Y(v_j)$
is relatively semiample. Since every $t_j>0$, its morphism is the Stein
factor of the product of the vertex morphisms. Therefore all $c$ in its relative interior
have the same ample model. Take $c=b$.
\end{proof}

\begin{corollary}\label{relative-model}
Under the hypotheses of Theorem~\ref{model-theorem}, the ample model of the
real adjoint exists over all of $S$ and restricts to the ample model on
every fibre. For rational $\mathcal B$ it is the analytic relative
$\operatorname{Proj}$ of the locally finitely generated algebra
$\mathcal R=\bigoplus_{m\geq0}p_*\cO_{\X}(mK_{\X/S}+\lfloor m\mathcal B\rfloor)$.
\end{corollary}
\begin{proof}
For each $s\in S$, by part (1) there is a good minimal model over a neighbourhood
$U$ of $s$. Let $\Z_U$ be the ample model of its semiample adjoint.
If $U,V$ are two such neighbourhoods, their simplices from the proof of
(1) have interiors containing $b$ in the same rational affine subspace.
We pick a rational $c$ in their intersection. The last paragraph of the
proof of (1) shows that $b$ and $c$ have the same ample model on
each neighbourhood. Hence $\Z_U$ and $\Z_V$ identify over $\X_{U\cap V}$
with the ample model of $\mathcal D(c)$. These identifications commute
with the maps from $\X$ and satisfy the cocycle condition, so the $\Z_U$
glue to a normal proper space $\Z$ over $S$.

Let $t\in U$ and pick a rational $c$ in the relative interior of the simplex for $U$. The common
subdivision in Lemma~\ref{adjoint-box} and the convexity argument in (1),
applied on $X_t$, show that $D_t(c)$ and $D_t(b)$ have the same ample
model. We denote by $\phi\colon\X_U\dashrightarrow\Y$ the model in the proof
of (1), let $p_\Y\colon\Y\to U$ be its structure map, and put
$\mathcal D_\Y(c)=K_{\Y/U}+\phi_*\mathcal B(c)$. We pick $q$
so that $q\mathcal B(c)$ is integral and $q\mathcal D_\Y(c)$ is free over
$U$. Negativity and the projection formula identify
$$
 \bigoplus_{n\geq0}p_*\cO_\X(nq\mathcal D(c))
 \quad\text{with}\quad
 \bigoplus_{n\geq0}(p_\Y)_*\cO_\Y(nq\mathcal D_\Y(c)).
$$
By Theorem~\ref{main} applied to $\mathcal B(c)$, base change holds in
every degree. So the relative $\operatorname{Proj}$ restricts to the
canonical model of $(X_t,B_t(c))$, and hence to the ample model for $b$.

If $\mathcal B$ is rational, we take $c=b$ above and let $q$ clear
its denominators. The displayed algebra is then the $q$-Veronese of
$\mathcal R$. The pieces of $\mathcal R$ are coherent, and over connected
open sets its section algebra is a domain. By \cite[Lemma~2.26 and
Remark~2.27]{Fuj}, finite generation of this Veronese implies that
$\mathcal R$ is locally finitely generated with the same relative
$\operatorname{Proj}$. This proves Theorem~\ref{model-theorem}(2).
\end{proof}

\begin{corollary}\label{gmm-compact-big}
Suppose in Theorem~\ref{main} that $S$ is compact and $\mathcal D$ is big
over $S$. Then $(\X,\mathcal B)$ has a good minimal model over $S$.
\end{corollary}
\begin{proof}
Let $Z\to S$ be the canonical model of Corollary~\ref{relative-model}.
Since $\mathcal D$ is big over $S$, the map $\X\dashrightarrow Z$ is bimeromorphic. Take a finite cover of $S$ by the neighbourhoods in the proof of Corollary~\ref{relative-model}. The relative interiors of the corresponding simplices contain a common neighbourhood of $b$ in its smallest rational affine subspace. Within this intersection choose a rational simplex containing $b$ in its relative interior. Each vertex has the same marked ample model $Z$ locally, and bimeromorphicity identifies its ample divisor with the transformed adjoint. Hence the vertex adjoints on $Z$ are
$\Q$-Cartier and ample over $S$. Hence $K_{Z/S}+B_Z$ is $\mathbb R$-Cartier
and ample over $S$, where $B_Z$ is the transform of $\mathcal B$. By the local
construction, $\X\dashrightarrow Z$ is $\mathcal D$-nonpositive. Among the
prime divisors $P$ contracted from $\X$, we keep those with
$$
 a(P,Z,B_Z)=a(P,\X,\mathcal B).
$$
By compactness there are finitely many of them. Resolving the inverse map
\cite[\S1.2.4]{AKMW}, we get a projective log resolution $r\colon W\to Z$ dominating $\X$, on which all the kept primes appear.
On $W$ we take the strict transform of $B_Z$, give each kept prime $P$
coefficient $-a(P,Z,B_Z)$, and give every other $r$-exceptional prime
coefficient $1-\varepsilon$. Since $(Z,B_Z)$ is klt, for a small rational
$\varepsilon>0$ we obtain a klt boundary $\Theta$ with
$$
 K_W+\Theta=r^*(K_Z+B_Z)+N,\qquad N\ge0,
$$
where $\Supp N$ consists precisely of the remaining $r$-exceptional
primes. Since $r$ is bimeromorphic, $K_W+\Theta$ is
$r$-big, and both spaces are compact. So by \cite[Theorem~2.15(3),(4)]{HLX},
running an MMP with scaling of an $r$-ample line bundle, we obtain a strongly
$\Q$-factorial klt model $g\colon Y\to Z$. On $Y$ the transform
$N_Y=K_Y+\Theta_Y-g^*(K_Z+B_Z)$ is effective, $g$-exceptional and
$g$-nef, hence zero by the negativity lemma. For a kept prime $P$,
$$
 a(P,W,\Theta)=a(P,Z,B_Z)=a(P,Y,\Theta_Y).
$$
By strict negativity of the MMP, $P$ is not contracted, while the primes in $\Supp N$ are contracted since $N_Y=0$. Hence $\Theta_Y$ is the transform of
$\mathcal B$, and every prime lost from $\X$ has strictly larger
discrepancy on $Y$. Finally, $K_{Y/S}+\Theta_Y=g^*(K_{Z/S}+B_Z)$
is relatively semiample.
\end{proof}

\appendix
\section{A family with no global good minimal model}\label{app-schoen}

\begin{proposition}\label{global-counterexample}
There is a smooth proper isotrivial family over an elliptic curve $S$ with
smooth projective threefold fibres and a boundary $\mathcal B=\tfrac12\mathcal E$,
where $\mathcal E$ is smooth over $S$. It has a relative good minimal model
over every sufficiently small disc. Over $S$ there is no
$(K_{\X/S}+\mathcal B)$-nonpositive bimeromorphic contraction to a normal
space on which the transformed adjoint is relatively semiample.
\end{proposition}
\begin{proof}
Let $f\colon R\to\mathbf P^1$ and $f'\colon R'\to\mathbf P^1$ be
rational elliptic surfaces from general cubic pencils with disjoint
discriminants. The Mordell--Weil rank of $R$ is $10-2=8$, and translation
by a nontorsion section $P$ extends to an automorphism $\tau$ of $R$
\cite[Corollary~6.13 and \S7.6]{SS}. For the zero section $O'$ of $R'$ put
$$
 V=R\times_{\mathbf P^1}R',\qquad q\colon V\to R',\qquad E=q^{-1}(O')\simeq R .
$$
The threefold $V$ is smooth with $K_V\sim0$ by adjunction \cite{Schoen}.
Since $\cO_E(E)=f^*\cO_{\mathbf P^1}(-1)$, restricting to $E$ we see that
\begin{equation}\label{schoen-fixed}
 H^0(V,jE)=\C\,\sigma_E^j \qquad (j\geq1).
\end{equation}

Let $S=\C/\Lambda$ be an elliptic curve and pick a surjection
$\chi\colon\Lambda\to\mathbb Z$. We descend $(V,E)\times\C$ by the action
$\gamma(v,z)=((\tau^{\chi(\gamma)},1)v,z+\gamma)$ to $(\X,\mathcal E)$.
Over a disc we take the product of a good minimal model of $(V,E/2)$
\cite{KMM,KMM04} and apply Lemma~\ref{gmm-factorial}. The projection $q$
descends to $Q\colon\X\to S\times R'$. Let $F$ be a general elliptic
fibre of $\mathcal E\to S\times O'$.
For $L\in\operatorname{Pic}(\mathcal E)$, monodromy fixes $c_1(L|_R)$.
Since $H^1(R,\cO_R)=0$, it fixes $L|_R$. On the generic elliptic curve
of $f$, translation by $P$ acts on $\operatorname{Pic}^d$ by
translation by $\pm dP$. Hence
$$
 d=\deg(L|_F)\quad\Longrightarrow\quad dP=0\quad\Longrightarrow\quad d=0 .
$$
In particular no effective divisor on $\mathcal E$ dominates $S\times O'$.
Now let $\phi\colon\X\dashrightarrow Y$ be a contraction as in the
statement. By \eqref{schoen-fixed}, the fixed part in degree $2j$ is
$j\mathcal E$ locally over $S$. Nonpositivity and semiampleness force
$\phi$ to contract $\mathcal E$ to a centre $C$ of dimension at most two.
Since the indeterminacy locus has codimension at least two in $\X$, it does not meet a general $F$: otherwise its intersection with $\mathcal E$ would
contain a divisor dominating $S\times O'$. So $\phi$ is holomorphic
near $F$.
If $\phi(F)$ is a curve, then $\dim C=2$. We resolve $C$ by $\widetilde C$
and the map $\mathcal E\dashrightarrow\widetilde C$ by $u\colon W\to\mathcal E$,
$v\colon W\to\widetilde C$. A general $F$ avoids the exceptional image
of $u$. Its image in $\widetilde C$ is a smooth fibre component dominated
by the rational surface $R$, hence is $\mathbf P^1$. Writing $e>0$ for
the degree on $F$, we see that
$$
 0=\deg\!\left((u_*v^*\omega_{\widetilde C/S})^{**}|_F\right)=-2e .
$$
Hence $\phi(F)$ is a point. For a Stein neighbourhood $U$ of $\phi(F)$, its inverse image in the holomorphic locus of $\phi$ contains $F$. By properness of $Q$, it contains $Q^{-1}(A)$ for some neighbourhood $A$ of $Q(F)$. The compact elliptic fibres over $A$ map to points of $U$. This contradicts bimeromorphicity.
\end{proof}

One can also take $\kappa=2$. We blow down $O'$ by $\beta\colon R'\to P_0$.
For a general very ample $A$ avoiding $\beta(O')$, the smooth divisor
$H=(\beta q)^*A$ is disjoint from $E$ and descends. Since $O'$ is
$\beta$-exceptional and $q$ has connected fibres,
$$
 H^0(V,j(E+H))=\sigma_E^j\,(\beta q)^*H^0(P_0,jA).
$$
So $(E+H)/2$ has the same fixed part, and its Iitaka morphism is
induced by $\beta q$ on the whole family.

\section*{Acknowledgments}

The first named author is supported by Overseas High-Level Young Talent Recruitment
Programs and is partially supported by the Chongqing Natural Science Foundation Innovation
and Development Joint Fund CSTB2023NSCQ-LZX0031. The second named author is partially
supported by the National Key Research and Development Program of China (Grant No.
2024YFA1013201), the Natural Science Foundation of Chongqing (Grant No.
CSTB2024NSCQ-LZX0040), and the Special Project of Chongqing Municipal Science and Technology
Bureau (Grant No.\ 2025CCZ015).


\begin{thebibliography}{ELMNP23}
\bibitem[AKMW02]{AKMW} D. Abramovich, K. Karu, K. Matsuki and J. W{\l}odarczyk,
\emph{Torification and factorization of birational maps}, J. Amer. Math. Soc.
\textbf{15} (2002), \href{https://doi.org/10.1090/S0894-0347-02-00396-X}{531--572}.
\bibitem[BP12]{BP} B. Berndtsson and M. P\u{a}un, \emph{Quantitative extensions
of pluricanonical forms and closed positive currents}, Nagoya Math. J.
\textbf{205} (2012), 25--65.
\bibitem[Bi83]{Bingener} J. Bingener, \emph{On deformations of K\"ahler spaces II},
Arch. Math. (Basel) \textbf{41} (1983), 517--530.
\bibitem[BM88]{BM} E. Bierstone and P. D. Milman,
\emph{Semianalytic and subanalytic sets}, Publ. Math. IH\'ES
\textbf{67} (1988), \href{https://doi.org/10.1007/BF02699126}{5--42}.
\bibitem[BCHM10]{BCHM} C. Birkar, P. Cascini, C. D. Hacon and J. McKernan,
\emph{Existence of minimal models for varieties of log general type},
J. Amer. Math. Soc. \textbf{23} (2010),
\href{https://doi.org/10.1090/S0894-0347-09-00649-3}{405--468}.
\bibitem[Bou02]{Bou} S. Boucksom, \emph{On the volume of a line bundle},
Internat. J. Math. \textbf{13} (2002), 1043--1063.
\bibitem[CP23]{CP} J. Cao and M. P\u{a}un, \emph{Infinitesimal extension of
pluricanonical forms}, \href{https://arxiv.org/abs/2012.05063v5}{arXiv:2012.05063v5} (2023).
\bibitem[Con06]{Con} B. Conrad, \emph{Relative ampleness in rigid geometry},
Ann. Inst. Fourier (Grenoble) \textbf{56} (2006),
\href{https://doi.org/10.5802/aif.2207}{1049--1126}.
\bibitem[DHY26]{DHY} O. Das, C. D. Hacon and J. I. Y\'a\~nez, \emph{MMP for
generalized pairs on K\"ahler $3$-folds}, Ann. Sc. Norm. Super. Pisa Cl. Sci. (5),
\href{https://doi.org/10.2422/2036-2145.202506_004}{published online (2026)}.
\bibitem[Dem10]{Dem} J.-P. Demailly,
\href{https://www-fourier.univ-grenoble-alpes.fr/~demailly/manuscripts/analmeth_book.pdf}{\emph{Analytic methods in algebraic geometry}},
Surveys of Modern Mathematics, vol.~1, Higher Education Press, 2010.
\bibitem[ELMNP23]{ELMNP} L. Ein, R. Lazarsfeld, M. Musta\c{t}\u{a},
M. Nakamaye and M. Popa, \emph{Erratum to the paper: Asymptotic invariants
of base loci}, \href{https://arxiv.org/abs/2309.16722}{arXiv:2309.16722} (2023).
\bibitem[FLR26]{FLR} J. Feng, Y. Li and S. Rao, \emph{Aspects of deformation
invariance of plurigenera of K\"ahler varieties},
\href{https://arxiv.org/abs/2610.17817v1}{arXiv:2610.17817v1}.
\bibitem[Fuj13]{Fuj13} O. Fujino, \emph{A transcendental approach to Koll\'ar's
injectivity theorem II}, J. Reine Angew. Math. \textbf{681} (2013),
\href{https://doi.org/10.1515/crelle-2012-0036}{149--174}.
\bibitem[Fuj15]{Fuj15} O. Fujino, \emph{Some remarks on the minimal model program
for log canonical pairs}, J. Math. Sci. Univ. Tokyo \textbf{22} (2015),
\href{https://www.ms.u-tokyo.ac.jp/journal/pdf/jms220105.pdf}{149--192}.
\bibitem[Fuj22]{Fuj} O. Fujino, \emph{Minimal model program for projective morphisms
between complex analytic spaces}, \href{https://arxiv.org/abs/2201.11315v1}{arXiv:2201.11315v1} (2022).
\bibitem[Fuj23]{Fuj23}
O. Fujino, \emph{Relative Bertini type theorem for multiplier
ideal sheaves}, Osaka J. Math. \textbf{60} (2023),
\href{https://doi.org/10.18910/89998}{207--226}.
\bibitem[HLR26]{HLR} C. D. Hacon, Y. Li and S. Rao, \emph{On pseudo-effectivity and
volumes of adjoint classes in K\"ahler families with projective central fiber},
to appear in Math. Ann., \href{https://arxiv.org/abs/2602.04158v3}{arXiv:2602.04158v3}.
\bibitem[HLX26]{HLX} C. D. Hacon, Y. Li and L. Xie,
\emph{Fujiki class $\mathcal C$ varieties and a K\"ahler criterion},
preprint (2026), \href{https://arxiv.org/abs/2608.20588v2}{arXiv:2608.20588v2}.
\bibitem[HMX13]{HMX13} C. D. Hacon, J. McKernan and C. Xu, \emph{On the
birational automorphisms of varieties of general type}, Ann. of Math. (2)
\textbf{177} (2013), 1077--1111.
\bibitem[HMX18]{HMX} C. D. Hacon, J. McKernan and C. Xu, \emph{Boundedness of moduli
of varieties of general type}, J. Eur. Math. Soc. \textbf{20} (2018), 865--901.
\bibitem[HW26]{HW} C. D. Hacon and P. Wang, \emph{Termination of the K\"ahler MMP
with scaling}, \href{https://arxiv.org/abs/2611.29683v1}{arXiv:2611.29683v1}.
\bibitem[HX26]{HX} C. D. Hacon and L. Xie, \emph{On the K\"ahler MMP and the
transcendental base-point-free theorem}, \href{https://arxiv.org/abs/2607.24986v2}{arXiv:2607.24986v2}.
\bibitem[Iit70]{Iit} S. Iitaka, \emph{Deformations of compact complex surfaces II},
J. Math. Soc. Japan \textbf{22} (1970), 247--261.
\bibitem[KMM94]{KMM} S. Keel, K. Matsuki and J. McKernan,
\emph{Log abundance theorem for threefolds}, Duke Math. J.
\textbf{75} (1994), 99--119.
\bibitem[KMM04]{KMM04} S. Keel, K. Matsuki and J. McKernan,
\emph{Corrections to \textquotedblleft Log abundance theorem for threefolds\textquotedblright},
Duke Math. J. \textbf{122} (2004),
\href{https://doi.org/10.1215/S0012-7094-04-12236-5}{625--630}.
\bibitem[KS60]{KS} K. Kodaira and D. C. Spencer, \emph{On deformations of complex
analytic structures. III. Stability theorems for complex structures},
Ann. of Math. (2) \textbf{71} (1960), 43--76.
\bibitem[Kol11]{Kol11} J. Koll\'ar, \emph{A local version of the
Kawamata--Viehweg vanishing theorem}, Pure Appl. Math. Q. \textbf{7} (2011),
\href{https://doi.org/10.4310/PAMQ.2011.v7.n4.a18}{1477--1494}.
\bibitem[Kol21]{Kol} J. Koll\'ar, \emph{Deformations of varieties of general type},
Milan J. Math. \textbf{89} (2021), 345--354.
\bibitem[KM92]{KM} J. Koll\'ar and S. Mori, \emph{Classification of three-dimensional
flips}, J. Amer. Math. Soc. \textbf{5} (1992), 533--703.
\bibitem[KNX18]{KNX} J. Koll\'ar, J. Nicaise and C. Xu,
\emph{Semi-stable extensions over 1-dimensional bases}, Acta Math. Sin.
(Engl. Ser.) \textbf{34} (2018), \href{https://doi.org/10.1007/s10114-017-7048-8}{103--113}.
\bibitem[Lev83]{Lev83} M. Levine, \emph{Pluri-canonical divisors on K\"ahler
manifolds}, Invent. Math. \textbf{74} (1983), 293--303.
\bibitem[Lev85]{Lev85} M. Levine, \emph{Pluri-canonical divisors on K\"ahler
manifolds. II}, Duke Math. J. \textbf{52} (1985), 61--65.
\bibitem[LPWZ26]{LPWZ} Y. Li, M. P\u{a}un, Z. Wang and X. Zhou, \emph{Psh
variation of singular K\"ahler--Einstein metrics and invariance of plurigenera for
smooth K\"ahler families}, \href{https://arxiv.org/abs/2610.07826v1}{arXiv:2610.07826v1}.
\bibitem[LW26]{LW} K. Liu and X. Wan, \emph{Invariance of logarithmic plurigenera
and adjoint volumes in K\"ahler families}, \href{https://arxiv.org/abs/2610.26519v1}{arXiv:2610.26519v1}.
\bibitem[Nak86]{Nak} N. Nakayama, \emph{Invariance of the plurigenera of
algebraic varieties under minimal model conjectures}, Topology \textbf{25}
(1986), 237--251.
\bibitem[RT22]{RT} S. Rao and I.-H. Tsai, \emph{Invariance of plurigenera and
Chow-type lemma}, Asian J. Math. \textbf{26} (2022), 507--554.
\bibitem[Sch88]{Schoen} C. Schoen, \emph{On fiber products of rational elliptic
surfaces with section}, Math. Z. \textbf{197} (1988), 177--199.
\bibitem[SS10]{SS} M. Sch\"utt and T. Shioda, \emph{Elliptic surfaces},
in \emph{Algebraic Geometry in East Asia---Seoul 2008}, Adv. Stud.
Pure Math. \textbf{60} (2010), 51--160.
\bibitem[Siu98]{Siu98} Y.-T. Siu, \emph{Invariance of plurigenera},
Invent. Math. \textbf{134} (1998), 661--673.
\bibitem[Siu02]{Siu} Y.-T. Siu, \emph{Extension of twisted pluricanonical sections
with plurisubharmonic weight and invariance of semipositively twisted plurigenera
for manifolds not necessarily of general type}, in \emph{Complex Geometry},
Springer, Berlin, 2002, 223--277.
\bibitem[Wei94]{Wei} C. A. Weibel, \emph{An introduction to homological algebra},
Cambridge Studies in Advanced Mathematics \textbf{38}, Cambridge University Press,
\href{https://doi.org/10.1017/CBO9781139644136}{1994}.
\bibitem[Zhu25]{Zhu} M. Zhu, \emph{Boundedness of stable minimal models with
klt singularities}, Int. Math. Res. Not. \textbf{2025}, no.~2,
\href{https://doi.org/10.1093/imrn/rnae293}{rnae293}.
\end{thebibliography}
\end{document}